\documentclass{rsepublic}
\usepackage{graphicx,subfigure}

\usepackage[english]{babel}
\usepackage[utf8]{inputenc}
\usepackage[T1]{fontenc}
\usepackage{amsmath}
\usepackage{amsfonts}
\usepackage{mathrsfs}
\usepackage{amssymb}
\usepackage{amsthm}
\usepackage{comment}
\usepackage{lmodern}

\usepackage[colorlinks=true]{hyperref}

\hypersetup{dvipdfmx,
       colorlinks,
       linkcolor=blue,
       filecolor=red,
       urlcolor=blue,
       citecolor=red,
       pdftitle={Concentration structures},
       pdfauthor={ Jun YANG},
       pdfsubject={   PDEs         },
       pdfkeywords={   Schrodinger equation},
       pdfproducer={ps2pdf}}
\numberwithin{equation}{section}
\numberwithin{table}{section}

{\theoremstyle{plain}

\usepackage{CJKnumb}

\newtheorem{remark}{\indent Remark}[section]

\numberwithin{equation}{section}

\newcommand{\R}{\mathbb{R}}

\date{} %\linespread{2.0}

\allowdisplaybreaks

\begin{document}

\title[Non-degeneracy of vortex solutions] {On the non-degeneracy of radial vortex solutions for a coupled Ginzburg-Landau system}
\author[L. Duan and J. Yang]
{\textbf{Lipeng Duan}%\thanks{dd}
\\
School of Mathematics and Statistics, Central China Normal University,
Wuhan 430079, P. R. China.
Email: ahudlp@sina.com
\\
\medskip
\textbf{Jun Yang}\footnote{Corresponding author}
\\
School of Mathematics and Statistics \& Hubei Key Laboratory of Mathematical
Sciences, Central China Normal University,
Wuhan 430079, P. R. China.
Email: jyang@mail.ccnu.edu.cn
}

%\author{Jun Yang}
%\address{Jun Yang,
%\newline\indent School of Mathematics and Statistics \& Hubei Key Laboratory of Mathematical
%\newline\indent Sciences, Central China Normal University,
%Wuhan 430079, P. R. China.
%}
%\email{jyang@mail.ccnu.edu.cn}

%\date{\today}

\MSdates{}
%\MSdates[¡®Received date¡¯]{¡®Accepted date¡¯}
\label{firstpage}

\maketitle

\begin{abstract}
For the following Ginzburg-Landau system in ${\mathbb R}^2$
 \begin{align*}
 \begin{cases}
 -\Delta w^+ +\Big[A_+\big(|w^+|^2-{t^+}^2\big)+B\big(|w^-|^2-{t^-}^2\big)\Big]w^+=0,
 \\[3mm]
-\Delta w^- +\Big[A_-\big(|w^-|^2-{t^-}^2\big)+B\big(|w^+|^2-{t^+}^2\big)\Big]w^-=0,
\end{cases}
\end{align*}
with constraints $ A_+, A_->0$, $B<0$, $B^2<A_+A_-$ and $t^+, t^->0$,
we will concern its linearized operator ${\mathcal L}$  around the radially symmetric solution $w(x)=(w^+, w^-): {\mathbb R}^2 \rightarrow\mathbb{C}^2$ of degree pair $(1, 1)$  and prove the non-degeneracy result: the kernel
of ${\mathcal L}$ is spanned by $\big\{\frac{\partial w}{\partial{x_1}}, \frac{\partial w}{\partial{x_2}}\big\}$ in a natural Hilbert space.
As an application of the non-degeneracy result, a solvability theory for the linearized operator~${\mathcal L}$ will be given.
\vspace{0.2cm}

%\noindent  {\bf 2000 MR Subject Classification:} \ 35Q35;
%76D05\\
\end{abstract}

%\footnote[0]
%{   \\
%\textit{}
% }
%\thanks{Corresponding author: J. Yang}

\section{Introduction}
\subsection{Existence of symmetric vortex solutions}
We consider the Ginzburg-Landau system in $\mathbb {R}^2$
\begin{align} \label{1.1}
 \begin{cases}
 -\Delta w^+ +\big[A_+(|w^+|^2-{t^+}^2)+B(|w^-|^2-{t^-}^2)\big]w^+=0,
 \\[2mm]
-\Delta w^- +\big[A_-(|w^-|^2-{t^-}^2)+B(|w^+|^2-{t^+}^2)\big]w^-=0,
\end{cases}
\end{align}
where $w=(w^+, w^-):\mathbb {R}^2\rightarrow \mathbb{C}^2$ is a complex vector-valued function. The Ginzburg-Landau system of this type
has been introduced in physical modes of Bose-Einstein  condensates and for the superconductors.
The reader can refer to the papers~\cite{AG}, \cite{KMM},~\cite{AB}
about the applications of the Ginzburg-Landau system.

\medskip
The energy functional of the system (\ref{1.1}) is
 \begin{align*}
{\mathbf E}( w)=\int_{\mathbb{R}^2}\frac{1}{2}|\nabla w|^2+&\int_{\mathbb{R}^2}\frac{1}{4}\Big[\,A_+(| w^+|^2-{t^+}^2)^2+A_-(| w^-|^2-{t^-}^2)^2
\\[1mm]
&\qquad\quad+2B(| w^+|^2-{t^+}^2)(| w^-|^2-{t^-}^2)\,\Big],
\end{align*}
where $ w=( w^+,  w^-):\mathbb {R}^2\rightarrow \mathbb{C}^2$.
For the balanced case, $A_+=A_-$ and $t^+=t^-=1/\sqrt{2}$, S. Alama, L. Bronsard and P. Minronescu proved the existence, uniqueness and monotonicity of radial
solutions in \cite{AlamaBronsardMironescu1}.
By making the following assumptions
\begin{align}
\label{H} A_+,A_->0,\qquad B^2<A_+A_-,\qquad t^+, t^->0,\tag{H1}
\end{align}
S. Alama and Q. Gao in \cite{AG}  concerned the general case and showed  that the energy functional of the  system (\ref{1.1}) is continuous, coercive, and convex, and then they proved that ${\mathbf E}(w)$ has a radially  symmetric minimizer~$w_n=(w_n^+, w_n^-)$  in the forms
\begin{align}\label{1.3}
w_n^+(r, \theta)=U_n^+(r)e^{in_+\theta},
\qquad
w_n^-(r, \theta)=U_n^-(r)e^{in_-\theta},
\end{align}
with given degree pair $n=(n_+,n_-)\in \mathbb{Z}^2$.
Here and in the sequel, we will use the convention
$$
x_1=r\cos\theta,
\qquad
x_2=r\sin\theta.
$$
By substituting (\ref{1.3}) in (\ref{1.1}), we can obtain an ODE system for the vortex profiles $U_n^{\pm}$ for $r\in (0, \infty)$,
\begin{align}\label{10.2}
 \begin{cases}
-{U_n^+}''-\frac{1}{r}{U_n^+}'+\frac{n_+^2}{r^2}U_n^++\Big[\,A_+({U_n^+}^2-{t^+}^2)+B({U_n^-}^2-{t^-}^2)\,\Big]U_n^+=0,
\vspace{3mm}
\\
-{U_n^-}{''}-\frac{1}{r}{U_n^-}'+\frac{n_-^2}{r^2}U_n^-+\Big[\,A_-({U_n^-}^2-{t^-}^2)+B({U_n^+}^2-{t^+}^2)\,\Big]U_n^-=0,
 \end{cases}
\end{align}
with the conditions
\begin{align}
U_n^\pm\geq 0  \quad \text{for~all } r\geq 0,
\qquad
U_n^\pm\rightarrow t^{\pm}\quad \mbox{as } r\rightarrow \infty,
\label{boundary1}
\\[2mm]
U_n^\pm(0)=0\quad \mbox{if } n_{\pm}\neq 0,
\qquad
{U_n^\pm}'(0)=0\quad \mbox{if } n_{\pm}= 0.
\label{boundary2}
\end{align}

\medskip
For the convenience of the readers, we give the following existence result and fundamental properties of solutions to (\ref{10.2})-(\ref{boundary2}).
The readers can refer to Proposition 1.1 in \cite{AG} for more details.

\begin{proposition}\cite{AG}\label{0016}
Let (\ref{H}) hold. Then there exists a unique solution $( U_n^+, U_n^-)$ to (\ref{10.2})-(\ref{boundary2}). Moreover, we have
$$
U_n^\pm\in C^\infty(0,\infty),
\qquad
U_n^\pm>0  ~\text{for~all }~r>0,
$$
$$
U_n^\pm\sim r^{|n_{\pm}|} ~\text{as }~r\rightarrow 0,
\qquad
\int^\infty_0 \big|{U_n^\pm}' \big|^2r\,\mathrm{d} r<\infty.
$$
In particular, $w_n=(U_n^+e^{in_+\theta}, U_n^-e^{in_-\theta})$ is an entire solution of (\ref{1.1}) in $\mathbb R^2$.
\qed
\end{proposition}

\medskip
The asymptotic behavior of vortex profiles $U_n^\pm$ as $r\rightarrow \infty$ are derived from the following proposition, see Theorem 3.2 in \cite{AG}.

\begin{proposition}\cite{AG}\label{0017}
Let (\ref{H}) hold. Suppose that $(U_n^+,U_n^-)$ is the solution of (\ref{10.2})-(\ref{boundary2}). Then we have
$$
U_n^\pm=t^{\pm}+\frac{\tilde c_\pm}{2r^2}+O(r^{-4})%~~~\text{as}~r\rightarrow \infty,$$
\quad\mbox{and}\quad
 {U_n^\pm}'=-\frac{\tilde c_{\pm}}{r^3}+O(r^{-5}),~~~\text{as}~r\rightarrow \infty,
$$
with
$$
\tilde c_{\pm}=\frac{B n_\mp^2-A_\mp n_\pm^2}{\big(A_+ A_--B^2\big)t^{\pm}}.
$$
\qed
\end{proposition}

\medskip
The monotonicity properties of the vortex profiles were also considered in Theorem 1.2 of \cite{AG}.
The validity of the properties is strongly dependent on the  interaction coefficient $B$.
\begin{proposition}\cite{AG}
Let (\ref{H}) hold. Assume that $w_n=(U_n^+e^{in_+\theta},U_n^-e^{in_-\theta})$ is the entire equivariant solution of (\ref{1.1}), which is described in Proposition \ref{0016}.
\medskip\\
{\textbf {(1).}}~If $B<0$, then ${U_n^\pm}'(r)\geq 0$ for all $r>0$ and any degree pair $n=(n_+,n_-)$.
\medskip\\
{\textbf {(2).}}~If $B>0,\, n_+\geq 1$ and $n_-=0$,  then ${U_n^+}'(r)\geq 0$ and ${U_n^-}'(r) \leq 0$ for all $r>0$.
\medskip\\
{\textbf {(3).}}~For any pair $(n_+,n_-)$ with $n_+\neq 0\neq n_-$, there exists $B_0>0$ such that  ${U_n^\pm}'(r)\geq 0$ for all $r>0$ and all $B$ with $0\leq B\leq B_0$.
\qed
\end{proposition}

%\medskip
%In the work \cite{ABG1}, the authors studied the  p-Ginzburg-Landau energy functional
%  \begin{align}\label{0014}
%  E_p(u)=\int_{\mathbb{R}^2}|\nabla u|^p+\frac{1}{2}(1-|u|^2)^2,
%  \end{align}
%  with $p>2$.
%  They considered the radially symmetric minimizer, i.e., $u_p=f_p(r)e^{i\theta}$, the degree one radially symmetric
%  solution of
%  \begin{align}
%  \frac{p}{2}\nabla\cdot(|\nabla u|^{p-2}\nabla u)+u(1-|u|^2)=0
%  \end{align}
%  Then, by substituting $u=f(r)e^{i\theta}$ into functional (\ref{0014}), we can obtain easily
%  $$E_p(u)=2\pi I_p(f),$$
%  with
%  \begin{align}\label{0015}
%  I_p(f)=\int_0^\infty\Big\{\Big(|f'|^2+\frac{f^2}{r^2}\Big)^{\frac{p}{2}}+\frac{1}{2}(1-f^2)^2\Big\}r {\mathrm d}r
%  \end{align}
%  In summary, we have the following  existence result and fundamental properties of the radially symmetric minimizer to the functional (\ref{0015}),
%  which we can infer from the literature \cite{ABG1}.
%  \begin{proposition}
%  There exist a unique  minimizer $f_p$ to the functional (\ref{0015}). Moreover, we have
%  $$f_p\in C^\infty(0,\infty),$$
%\qquad
%$$0<f_p<1 ~\text{for~all }~r>0,$$
%\qquad
%$$f_p\sim r ~\text{as }~r\rightarrow 0,$$
%\qquad
%$${f_p}'>0~~\text{in }~(0,\infty),$$
%\qquad
%$$1-{f_p}^2\sim \frac{p}{2}\frac{1}{r^p}~\text{and }~{f_p}'\sim\frac{p^2}{4}\frac{1}{r^{p+1}}~\text{as }~r\rightarrow \infty.$$
%  \end{proposition}

\medskip
\subsection{Stability and non-degeneracy of symmetric vortex solutions}
The  stability or non-degeneracy of the minimizer  is of fundamental importance in analysis.
There are many works involving  non-degeneracy and  stability,  see, e.g., \cite{AG2}, \cite{ABG1},
\cite{CM}, \cite{Lin}, \cite{MP}, \cite{MS} and so on.

\medskip
For the classical (single complex component) Ginzburg-Landau equation on the disc $B_R(0)\subset{\R}^2$,
\begin{align}\label{001}
-\Delta u=u(1-|u|^2),
\end{align}
P. Mironescu \cite{MP}  considered  the radial solution of degree one and gave the stability result
in the sense that corresponding quadratic form is positive definite.
By using the Fourier decomposition method, T. C. Lin \cite{Lin} studied the  (single complex component) Ginzburg-Landau equation in $B_1(0)$ with
boundary condition
 \begin{align*}
-\Delta u=\frac{1}{\epsilon^2}u(1-|u|^2)\quad\text{in } B_1(0),
\qquad
u=g\quad \text{on } \partial B_1(0),
\end{align*}
and proved the  stability result for the radial symmetric solutions of degree one.
M. del Pino,  P. Felmer and M. Kowalczyk \cite{DFK} studied the classical (single complex component) Ginzburg-Landau equation  (\ref{001}) in $\mathbb{R}^2$, and they proved the non-degeneracy result for the radial symmetric vortex solution of degree one in a given Hilbert space.
For the solutions of degree $d\geq 2$ of  Ginzburg-Landau mode, there are stability or instability results (see,~e.g., \cite{CM}, \cite{MP}, \cite{MS}).

\medskip
 In \cite{ABG0} and \cite{ABG1}, Y. Almog, L. Berlyand, D. Golovaty and I. Shafrir considered the minimization of a $p$-Ginzburg-Landau energy functional
  \begin{align}\label{0018}
 {\mathbf E}_p(u)=\int_{\mathbb{R}^2}|\nabla u|^p+\frac{1}{2}(1-|u|^2)^2,
  \end{align}
  with $p>2$,
  over the class of radially symmetric functions of degree one. The radially symmetric minimizer, i.e., $u_p=f_p(r)e^{i\theta}$, is the
   radially symmetric solution of degree one to
  \begin{align}
  \frac{p}{2}\nabla\cdot(|\nabla u|^{p-2}\nabla u)+u(1-|u|^2)=0.
  \end{align}
 Moreover they studied   the existence, uniqueness and the asymptotic  behavior results for the radially symmetric minimizer.
 In particular, when $2<p\leq 4$,   the radially symmetric solution is locally stable in the sense that the second variation of the functional $E_p(u)$ at the
 minimizer $u_p=f_p(r)e^{i\theta}$ is positive and the kernel
of linearized operator ${\mathcal L}_p$ around the minimizer is spanned by $\{\frac{\partial u_p}{\partial{x_1}},\frac{\partial u_p}{\partial{x_2}},\frac{\partial u_p}{\partial{\theta}}\}$ in a natural Hilbert space.

\medskip
Let us go back to the system (\ref{1.1}).
For the balanced case, $A_+=A_-$ and $t^+=t^-=1/\sqrt{2}$, S. Alama, L. Bronsard and P. Minronescu concerned the stability and bifurcation of radial
solutions on a disc \cite{AlamaBronsardMironescu2}.
Recently,   S. Alama and Q. Gao \cite{AG2} considered the  Ginzburg-Landau system (\ref{1.1}) on the disk in $\mathbb{R}^2$ with a symmetric, degree-one boundary Dirichlet condition, and studied its stability for $B<0$ case and instability for case of positive $B$ close to zero in sense of the spectrum of the second variation of the energy.
By suitable rescaling, they formally proposed that the translation invariance of the entire radial solution will lead to the whole null space for the
second variation of energy if $B<0$.

\medskip
Whence, in the present paper, we will concern the non-degeneracy of the radial vortex solution with degree pair $(1, 1)$ for  (\ref{1.1})
and give an affirmative answer to the above problem.
%we focus on the special case of  Ginzburg-Landau system (\ref{1.1}) in the form
%\begin{align} \label{equation1.1}
% \begin{cases}
% -\Delta w^+ +\big[A_+(|w^+|^2-{t^+}^2)+B(|w^-|^2-{t^-}^2)\big]w^+=0,
% \\[2mm]
%-\Delta w^- +\big[A_-(|w^-|^2-{t^-}^2)+B(|w^+|^2-{t^+}^2)\big]w^-=0.
%\end{cases}
%\end{align}
In the rest part of this paper, we use the notation
\begin{align}\label{(w^+,w^-)}
w=(w^+,w^-)
\quad\mbox{with}\quad
w^\pm=U^\pm e^{i\theta},
\end{align}
to denote the radially symmetric vortex solution  with degree pair $(n_+,n_-)=(1,1)$.
Then we get the corresponding ODE system
\begin{align}\label{1.2}
 \begin{cases}
-{U^+}''-\frac{1}{r}{U^+}'+\frac{1}{r^2}{U^+}+\Big[\,A_+({U^+}^2-{t^+}^2)+B({U^-}^2-{t^-}^2)\,\Big]U^+=0,\vspace{3mm}
\\
-{U^-}''-\frac{1}{r}{U^-}'+\frac{1}{r^2}{U^-}+\Big[\,A_-({U^-}^2-{t^-}^2)+B({U^+}^2-{t^+}^2)\,\Big]U^-=0.
 \end{cases}
\end{align}
Together with  Proposition \ref{0016} and Proposition \ref{0017}, when $(n_+,n_-)=(1,1)$,
the properties of the vortex profiles  $U^\pm$  are provided in the following
\begin{align}\label{0011}
\begin{cases}
&U^\pm\in C^\infty(0,\infty),
\qquad
0< U^\pm<t^\pm  ~\text{for~all }~r>0,
\\[2mm]
&U^\pm\sim r ~\text{as }~r\rightarrow 0,
\qquad
U^\pm\sim t^\pm-\frac{c_\pm}{2r^2}~\text{as }~r\rightarrow \infty,
\\[2mm]
&{U^\pm}'>0~\text{when}~B<0,
\qquad
{U^\pm}'\sim \frac{c}{r^3}~\text{as }~r\rightarrow \infty,
\end{cases}
\end{align}
where $c_\pm=\frac{A_\mp-B}{(A_+A_--B^2)t^\pm}$.
The linearized operator ${\mathcal L}=({\mathcal L}_+, {\mathcal L}_-)$  of (\ref{1.1}) around $w$ has the components
\begin{align}\label{0.1}
{\mathcal L}_\pm(\phi)=&-\Delta \phi^\pm \,+\,\Big[\,A_\pm\big({U^\pm}^2-{t^\pm}^2\big)+B\big({U^\mp}^2-{t^\mp}^2\big)\,\Big]\phi^\pm\nonumber
\\
&\,+\,2A_\pm {\mathbf {Re}}\Big(w^\pm\overline{\phi^\pm}\Big)w^\pm\,+\,2B{\mathbf {Re}}\Big(w^\mp \overline{\phi^\mp}\Big)w^\pm.
\end{align}
Direct computation shows that
\begin{align}\label{0.9}
{\mathcal L}\Big(\frac{\partial w}{\partial{x_1}}\Big)={\mathcal L}\Big(\frac{\partial w}{\partial{x_2}}\Big)={\mathcal L}(iw)=0.
\end{align}
The quadratic form given by the second variation of the energy functional $E$ around $w$ is
\begin{align}\label{0.2}
\boldsymbol{B}(\phi,\phi)=
&\int_{\mathbb{R}^2}|\nabla\phi|^2
\,+\,
\int_{\mathbb{R}^2}\Bigg[\,A_+\big({U^+}^2-{t^+}^2\big)+B\big({U^-}^2-{t^-}^2\big)\,\Bigg]|\phi^+|^2\nonumber
\\[2mm]
&\,+\,
\int_{\mathbb{R}^2}\Bigg[\,A_-\big({U^-}^2-{t^-}^2\big)+B\big({U^+}^2-{t^+}^2\big)\,\Bigg]|\phi^-|^2\nonumber
\\[2mm]
&+\int_{\mathbb{R}^2}\Bigg[\,2A_+\Big|{\mathbf {Re}}\big(\,\overline{w^+}\phi^+\big)\Big|^2+2A_-\Big|{\mathbf {Re}}\big(\,\overline{w^-}\phi^-\big)\Big|^2\,\Bigg]
\nonumber
\\[2mm]
&\,+\,
\int_{\mathbb{R}^2}4B {\mathbf {Re}}\big(\phi^+\overline{w^+}\,\big){\mathbf {Re}}\big(\phi^-\overline{w^-}\,\big).
\end{align}
By defining the inner product, for any $u=(u^+, u^-)$, $v=(v^+, v^-) : \mathbb{R}^2\rightarrow \mathbb{C}^2$,
\begin{align}\label{0.3}
\langle u,v\rangle=\langle u^+, v^+\rangle_{\mathcal R}\,+\,\langle u^-, v^-\rangle_{\mathcal R}
\quad\text{with}\quad
\langle u^\pm,v^\pm\rangle_{\mathcal R}={\mathbf {Re}}\int_{\mathbb{R}^2}u^{\pm}\overline{v^\pm},
\end{align}
 we can get
\begin{align*}
\boldsymbol{B}(\phi, \phi)=\langle {\mathcal L}(\phi), \phi\rangle=\langle {\mathcal L}_+(\phi),\phi^+\rangle_{\mathcal R}\,+\,\langle {\mathcal L}_-(\phi), \phi^-\rangle_{\mathcal R}.
\end{align*}
As a result of (\ref{0.9}), we  have
\begin{align}\label{Bnull}
\boldsymbol{B}(\phi,\phi)=0,
\end{align}
 for the special functions
\begin{align}\label{Bnull1}
\phi=\frac{\partial w}{\partial{x_1}}=\Big(\frac{\partial w^+}{\partial{x_1}}, \frac{\partial w^-}{\partial{x_1}}\Big)
 \quad\text{or}\quad
 \phi=\frac{\partial w}{\partial{x_2}}=\Big(\frac{\partial w^+}{\partial{x_2}}, \frac{\partial w^-}{\partial{x_2}}\Big),
\end{align}
where
\begin{align}
 \frac{\partial w^\pm}{\partial{x_1}}=\frac{1}{2}\Big({U^\pm}'-\frac{U^\pm}{r}\Big)e^{2i\theta}\,+\,\frac{1}{2}\Big({U^\pm}'+\frac{U^\pm}{r}\Big),
 \\[1mm]
 \frac{\partial w^\pm}{\partial{x_2}}=\frac{-i}{2}\Big({U^\pm}'-\frac{U^\pm}{r}\Big)e^{2i\theta}\,+\,\frac{i}{2}\Big({U^\pm}'+\frac{U^\pm}{r}\Big).
\end{align}
For the non-degeneracy, we mean the positivity of the quadratic form $\boldsymbol{B}(\phi,\phi)$ and that the kernel of the linearized operator ${\mathcal L}$ at $w$ is a linear combination of $\{\frac{\partial w}{\partial x_1},\frac{\partial w}{\partial x_2}\}$ in the space $\mathcal H$.
Because of (\ref{0.2}) and analogously to \cite{DFK}, we consider perturbations in a "natural" Hilbert space $\mathcal H$ with the norm as follows
\begin{align}
\begin{aligned}\label{mathcalH}
\|\phi\|^2_{\mathcal H}=
&\int_{\mathbb{R}^2}|\nabla\phi|^2
\,+\,\int_{\mathbb{R}^2}\Big[\,A_+\big({t^+}^2-{U^+}^2\big)-B\big({t^-}^2-{U^-}^2\big)\,\Big]|\phi^+|^2
\\
&\,+\,
\int_{\mathbb{R}^2}\Big[\,A_-\big({t^-}^2-{U^-}^2\big)-B\big({t^+}^2-{U^+}^2\big)\,\Big]|\phi^-|^2,
\end{aligned}
\end{align}
for any $\phi=(\phi^+, \phi^-): {\mathbb R}^2\rightarrow {\mathbb C}^2$.
In fact, as   $A_\pm > 0, B < 0,~0<U^\pm < t^\pm   $,  we can easily show  that  $\| \cdot\|$  is a norm.
Together with the  properties (\ref{0011}) of the vortex profiles $U^{\pm}$~of the minimizer $w$, we can know that $iw$ is not contained in the space $\mathcal H$.

\medskip
In the paper \cite{DFK},~the authors proved the  non-degeneracy of the solution of degree one for the classical (single complex component) Ginzburg-Landau equation in $\mathbb{R}^2$ by using Fourier decomposition method and some ODE techniques.
However, it is more complicated to get the  non-degeneracy for  the minimizer of the Ginzburg-Landau system (\ref{1.1}) due to the existence of the interaction terms by using the method in \cite{DFK}.
In this paper, we will overcome the difficulty by the method from \cite{ABG1} (using Picone's identities)
%Analogous to the arguments in \cite{AG2} and \cite{AB}, we make
under one more assumption
\begin{align}\label{H0}
B<0, \tag{H2}
\end{align}
due to technical reasons,
see the \textbf{Remark \ref{02}}.

\medskip
Now we give the main results.
\begin{theorem}\label{theorem1.1}
Assume that (\ref{H}) and (\ref{H0}) hold.
Suppose ${\mathcal L}(\phi)=0$ for some $\phi\in \mathcal H$.
Then we have
$$
\phi=c_1\frac{\partial w}{\partial x_1}+c_2\frac{\partial w}{\partial x_2},
$$
for certain constants $c_1$ and $c_2$.
\qed
\end{theorem}

\medskip
Note that we can use the same method to get the nondegeneracy results under the validity of (\ref{H}) and (\ref{H0})
if $w$ has degree pair $(-1, -1)$, $(-1, 1)$ or $(1, -1)$.
An interesting problem is to consider the validity of the non-degeneracy result for the minimizer of
Ginzburg-Landau system (\ref{1.1}) in $\mathbb{R}^2$ for $B>0$ case.
On the other hand, to the best knowledge of the authors, there are no non-degeneracy results for the vortex solutions with higher degrees
of the  Ginzburg-Landau equation with single complex component or a system like (\ref{1.1}).

\medskip
As an application of Theorem \ref{theorem1.1}, we give the following Fredholm alternative theorem for the linearized operator ${\mathcal L}$.

\begin{theorem}\label{theorem1.2}
Assume that (\ref{H}) and (\ref{H0}) hold.
Consider the equation
\begin{align}
\label{1.22} {\mathcal L}(\psi)=h~~~\text{in} ~\mathbb R^2,
\end{align}
with given $ h=(h^+,h^-) : {\mathbb R}^2\rightarrow {\mathbb C}^2$ satisfying
$$
\int_{\mathbb R^2}|h|^2(1+r^{2+\sigma})< +\infty
\quad\mbox{ for some } \sigma>0.
$$
If
\begin{align}
\label{1.35}\langle h^\pm,iw^\pm \rangle=\Big\langle h^\pm,\frac{\partial w^\pm}{\partial {x_1}} \Big\rangle=\Big\langle h^\pm,\frac{\partial w^\pm}{\partial {x_2}} \Big\rangle=0,
\end{align}
then (\ref{1.22}) has a solution $\psi\in \mathcal H$ which satisfies
\begin{align}\label{03}
\|\psi\|^2_{\mathcal H}\leq \int_{\mathbb R^2}|h|^2(1+r^{2+\sigma}).
\end{align}
Moreover, for any solution $\hat\psi^*\in \mathcal H$ of the equation (\ref{1.22}), $\psi$ has the form
\begin{align}\label{1.56}
\psi=\hat\psi^*+c_1\frac{\partial w}{\partial x_1}+c_2\frac{\partial w}{\partial x_2},
\end{align}
where $c_1$ and $c_2$ are two constants.
\qed
\end{theorem}

By following the work of \cite{DFK}, in Section \ref{section3} we will prove the solvability theory in Theorem \ref{theorem1.2}.
A solvability theory for the linearized operator is of crucial importance in the use of singular perturbation  methods for the construction of various solutions with much more complicated vortex structures of problems where the rescaled vortex $w$ provide a canonical profile.
More details about the subject for the single complex component Ginzburg-Landau equation are shown in the literature \cite{LinWeiYang},
\cite{PR},  \cite{WeiYang1} and \cite{Yang}.

\vspace{3mm}
\noindent{\textbf{Plan of the paper}}\\
We organize the paper as follows. In section 2, we will give the proof of  non-degeneracy  results by  using the method in the paper \cite{ABG1}.
In section 3, similar to the paper \cite{DFK}, based on  the non-degeneracy result, we will give a proof of the solvability theory~(Theorem \ref{theorem1.2}) for the  linear equation (\ref{1.22}).

\vspace{3mm}
\noindent{\textbf{Notation}}\vspace{1mm}\\
Throughout this paper,  we employ  $C, c ~\text{or}~ C_j,c_j,  j=0,1,2,\cdots$ to denote certain constants.
Furthermore,  we may make the abuse of notation by writing $\psi_j,\, \phi_j, j=0,1,2,\cdots$, from line to line in the present paper.

\medskip
\section{ Non-degeneracy  Results }\label{section2}

By following the method in \cite{ABG1},~now we give the proof of  non-degeneracy results.

\subsection{Preliminaries  for the proof of Theorem \ref{theorem1.1}}
The following decomposition by the Fourier series plays an important role in our study.
Now for any $\phi\in \mathcal{H}$ (see (\ref{mathcalH})),  we can decompose $\phi$ in its Fourier modes and write the quadratic form $\boldsymbol{B}(\phi,\phi)$ as a direct sum.
For any $\phi=(\phi^+,\phi^-): {\mathbb R}^2\rightarrow {\mathbb C}^2$, we write
\begin{align}\label{decompositionphi00}
\phi^{\pm}=\sum_{k\in \mathbb{Z}}\phi_k^{\pm}(r)e^{i k\theta},
\quad\mbox{or}\quad
\phi=\sum_{k\in \mathbb{Z}}\phi_k(r) e^{i k\theta}\ \mbox{with } \phi_k(r)=\big(\phi_k^+(r),\phi_k^-(r)\big),
\end{align}
where for any $k\in {\mathbb Z}$
\begin{align}\label{mathcalS}
\phi_k\in \mathcal S:=\Big\{\psi\in H_{loc}^1(\mathbb R_+, \mathbb C^2)\cap L_r^2(\mathbb R_+, \mathbb C^2)\,:\, \int_0^\infty\Big[\,|\psi'|^2+\frac{1}{r^2}|\psi|^2\,\Big]r\,{\mathrm d}r<+\infty\Big\}.
\end{align}
And then we get
\begin{align}\label{1.45}
\frac{1}{2\pi}\boldsymbol{B}(\phi,\phi)= 2\mathbb B_0 (\phi_1,\phi_1)+ \sum_{k=1}^{\infty}{\mathbb B}_k (\phi_{1+k},\phi_{1-k}),
 \end{align}
where
\begin{align}\label{1.46}
{\mathbb B}_0 (\phi_1,\phi_1)=&\int_0^\infty\Big[\,|(\phi_1^+)'|^2+\frac{1}{r^2}|\phi_1^+|^2\,\Big]r\,{\mathrm d}r
\,+\,\int_0^\infty \Big[\,|(\phi_1^-)'|^2+\frac{1}{r^2}|\phi_1^-|^2\,\Big]r\,{\mathrm d}r
\nonumber
\\[1mm]
&\,+\,\int_0^\infty 2\Big[\,A_+{U^+}^2|{\mathbf {Re}}(\phi_1^+)|^2+A_-{U^-}^2|{\mathbf {Re}}(\phi_1^-)|^2\,\Big]r\,{\mathrm d}r
\nonumber
\\[1mm]
&\,+\,\int_0^\infty 4B U^+U^-{\mathbf {Re}}(\phi_1^+){\mathbf {Re}}(\phi_1^-)r\,{\mathrm d}r
\nonumber\\[1mm]
&\,+\,\int_0^\infty\Big[\,A_+\big({U^+}^2-{t^+}^2\big)+B\big({U^-}^2-{t^-}^2\big)\,\Big]|\phi_1^+|^2r\,{\mathrm d}r
\nonumber\\[1mm]
&\,+\,\int_0^\infty\Big[\,A_-\big({U^-}^2-{t^-}^2\big)+B\big({U^+}^2-{t^+}^2\big)\,\Big]|\phi_1^-|^2r\,{\mathrm d}r,
\end{align}
and for $k\geq 1$
\begin{align}\label{1.20}
&{\mathbb B}_k (\phi_{1+k},\phi_{1-k})
\nonumber\\[1mm]
&=\int_0^\infty\Big[\,|(\phi_{1+k}^+)'|^2+|(\phi_{1-k}^+)'|^2+\frac{(1+k)^2}{r^2}|\phi_{1+k}^+|^2+\frac{(1-k)^2}{r^2}|\phi_{1-k}^+|^2\,\Big]r\,{\mathrm d}r
\nonumber\\[1mm]
&\quad\,+\,\int_0^\infty\Big[\,|(\phi_{1+k}^-)'|^2+|(\phi_{1-k}^-)'|^2+\frac{(1+k)^2}{r^2}|\phi_{1+k}^-|^2+\frac{(1-k)^2}{r^2}|\phi_{1-k}^-|^2\,\Big]r\,{\mathrm d}r
\nonumber\\[1mm]
&\quad\,+\,\int_0^\infty\Big[\,A_+{U^+}^2\big|\phi_{1+k}^++\overline {\phi_{1-k}^+}\big|^2+A_-{U^-}^2\big|\phi_{1+k}^-+\overline {\phi_{1-k}^-}\big|^2\,\Big]r\,{\mathrm d}r
\nonumber\\[1mm]
&\quad\,+\,\int_0^\infty 2B U^+U^-\Bigg[{\mathbf {Re}}\Big(\big(\phi_{1+k}^++\overline {\phi_{1-k}^+}\,\big)\overline {\phi_{1+k}^-}\Big)
+{\mathbf {Re}}\Big(\big(\phi_{1-k}^++\overline {\phi_{1+k}^+}\,\big) \overline{\phi_{1-k}^-}\,\Big)\Bigg]r\,{\mathrm d}r
\nonumber\\[1mm]
&\quad\,+\,\int_0^\infty\Big[\,A_+\big({U^+}^2-{t^+}^2\big)+B\big({U^-}^2-{t^-}^2\big)\,\Big](|\phi_{1+k}^+|^2+|\phi_{1-k}^+|^2)r\,{\mathrm d}r
\nonumber\\[1mm]
&\quad\,+\,\int_0^\infty\Big[\,A_-\big({U^-}^2-{t^-}^2\big)+B\big({U^+}^2-{t^+}^2\big)\,\Big](|\phi_{1+k}^-|^2+|\phi_{1-k}^-|^2)r\,{\mathrm d}r.
\end{align}
For the  proof of the above results,  the reader can refer  the paper \cite{AG2}, and we omit it here.

\medskip
The decompositions in (\ref{decompositionphi00}) are naturally associated to the functions
 $ \frac{\partial w}{\partial x_1}, \frac{\partial w}{\partial x_2}, iw$ in the kernel  of the linearized operator ${\mathcal L}$.
 In fact, we can write
\begin{align*}
iw=i\phi_1e^{i\theta},
\qquad
\frac{\partial w}{\partial{x_1}}=\phi_2e^{2i\theta}+\phi_0,
\qquad
\frac{\partial w}{\partial{x_2}}=-i\phi_2e^{2i\theta}+i\phi_0,
\end{align*}
where the functions have the special forms
\begin{align}
&\phi_1=({\phi_1}^+, {\phi_1}^-)\quad\text{with}\quad{\phi_1}^\pm=U^\pm,
\label{9}
\\[1mm]
&\phi_2=({\phi_2}^+, {\phi_2}^-)\quad\text{with}\quad{\phi_2}^\pm=\frac{1}{2}\Big({U^\pm}'-\frac{U^\pm}{r}\Big),
\label{10}
\\[1mm]
&\phi_0=({\phi_0}^+, {\phi_0}^-)\quad\text{with}\quad{\phi_0}^\pm=\frac{1}{2}\Big({U^\pm}'+\frac{U^\pm}{r}\Big).
\label{11}
\end{align}
%Note that \eqref{1.45} and \eqref{Bnull}-\eqref{Bnull1} imply that
%\begin{align}\label{mathbbBnull}
% \mathbb B_1 ( -\Phi_2, \Phi_0)=\frac{1}{\pi}\boldsymbol{B}\Big(\frac{\partial w}{\partial{x_1}}, \frac{\partial w}{\partial{x_1}}\Big)
%=0,
%\qquad
%\mathbb B_1 ( i\Phi_2, i\Phi_0)=\frac{1}{\pi}\boldsymbol{B}\Big(\frac{\partial w}{\partial{x_2}}, \frac{\partial w}{\partial{x_2}}\Big)
%=0.
%\end{align}

\medskip
As a result of (\ref{1.45}), we will discuss the quadratic forms $\{\mathbb B_k\}$ on $\mathcal S\times\mathcal S$ to get the non-degeneracy result.
It is easy to show that
$$
{\mathbb B}_0 (\phi_{1},\phi_{1})\geq 0
\quad\mbox{and}\quad
{\mathbb B}_k (\phi_{1+k},\phi_{1-k})>0,\ \forall\, k\geq 2,
$$
see Propositions \ref{proposition2.1} and \ref{proposition2.2}.
The key step is to establish that ${\mathbb B}_1 (\phi_{2},\phi_{0})\geq 0$,
and the equality holds if and only if $\phi_{2}, \phi_{0}$ have the special forms in (\ref{10})-(\ref{11}),
 see Proposition \ref{proposition2.3}.

\vspace{3mm}
\subsection{ Proof of Theorem \ref{theorem1.1}}
Before giving the complete proof of Theorem \ref{theorem1.1},
we analyze the quadratic forms  $\{\mathbb B_k\}$ separately and give several propositions.
%\noindent{\textbf{ The case $k=0$}}.

\begin{proposition}\label{proposition2.1}
We have that
$$
\mathbb B_0 (\phi, \phi) \geq  0,
\quad\forall\, \phi\in \mathcal S,
$$
and the equality holds if and only if $\phi= i (c_3\,U^+, c_4\,U^-)$  where $c_3, c_4$ are  real constants.
\end{proposition}

{\textit{Proof of Proposition \ref{proposition2.1}}}:
We can easily get that
$$
\mathbb B_0 (i|\phi|,i|\phi|)\leq \mathbb B_0 (\phi,\phi),
$$
for any $\phi\in \mathcal S $, and the strict inequality holds unless $\phi$ takes only purely imaginary  values.
Hence, it is natural to consider the  functional on ${\mathcal S} \times {\mathcal S}$
\begin{align*}
\tilde {\mathbb B}_0 (\phi,\phi)=&\int_0^\infty\Big[\,|(\phi^+)'|^2+\frac{1}{r^2}|\phi^+|^2\,\Big]r\,{\mathrm d}r
\,+\,\int_0^\infty \Big[\,|(\phi^-)'|^2+\frac{1}{r^2}|\phi^-|^2\,\Big]r\,{\mathrm d}r
\\[1mm]
&\,+\,\int_0^\infty\Big[\,A_+({U^+}^2-{t^+}^2)+B({U^-}^2-{t^-}^2)\,\Big]|\phi^+|^2r\,{\mathrm d}r
\\[1mm]
&\,+\,\int_0^\infty\Big[\,A_-({U^-}^2-{t^-}^2)+B({U^+}^2-{t^+}^2)\,\Big]|\phi^-|^2r\,{\mathrm d}r.
\end{align*}
For any $\phi\in C_c^\infty(0,\infty;\mathbb{R}^2)$,  whose support does not contain the origin, we write $\phi=(U^+\xi_1,U^-\xi_2)$  for some real functions $\xi_1, \xi_2$.
Integration by part together with (\ref{1.2})
yields
\begin{align}\label{1.19}
\tilde {\mathbb B}_0 (\phi,\phi)=\int_0^\infty ({U^+}^2{\xi_1'}^2+{U^-}^2{\xi_2'}^2)r\,{\mathrm d}r.
\end{align}

\medskip
The formula (\ref{1.19}) also holds for any smooth function $\phi$ whose support contains the origin.
The proof relies on a cutoff function. Using a density argument, (\ref{1.19}) is also true for any $\phi\in \mathcal H$. Then  if $\tilde {\mathbb B}_0 (\phi,\phi)=0$,  it is clearly that $\xi_1,\xi_2$ are  constants. This finishes the proof of Proposition \ref{proposition2.1}.
\qed

%\vspace{3mm}
%\noindent{\textbf{ The case $k\geq2$}}.
\medskip
\begin{proposition}\label{proposition2.2}
For each $k\geq 2$,   we have
$$ \mathbb {B}_k(\phi_1,\phi_2)>0 ~~~\text{for~all} ~(\phi_1,\phi_2)\in {\mathcal S} \times{\mathcal S} \setminus\{0,0\}. $$
\end{proposition}

{\textit{Proof of Proposition \ref{proposition2.2}}}:
From the definition of the functional  $\mathbb B_k(\phi_1,\phi_2) $ in (\ref{1.20}), we can know easily that
$$
\mathbb B_k(\phi_1,\phi_2)  \geq  \tilde {\mathbb B}_0(|\phi_1|,|\phi_1|)+\tilde {\mathbb B}_0(|\phi_2|,|\phi_2|)~~~\text{for~all} ~(\phi_1,\phi_2)\in {\mathcal S} \times{\mathcal S} \setminus\{0,0\}.
$$
By using Proposition \ref{proposition2.1}, we conclude that $\mathbb {B}_k(\phi_1,\phi_2)\geq0$ and equality hold if and only if $\phi_1=\phi_2=0$.
This finishes the proof of Proposition \ref{proposition2.2}.
\qed

\medskip
\medskip
%\noindent{\textbf{ The case $k=1$}}.
Now, we discuss the quadratic form  $\mathbb B_1(\phi_2,\phi_0)$, which is the most delicate one.
Note that
\begin{align}
\mathbb B_1 (\phi_{2},\phi_{0})=&\int_0^\infty\Big[\,|(\phi_{2}^+)'|^2+|(\phi_{0}^+)'|^2+\frac{4}{r^2}|\phi_{2}^+|^2\,\Big]r\,{\mathrm d}r
+\int_0^\infty \Big[\,|(\phi_{2}^-)'|^2+|(\phi_{0}^-)'|^2+\frac{4}{r^2}|\phi_{2}^-|^2\,\Big]r\,{\mathrm d}r
\nonumber\\[1mm]
&+\int_0^\infty \Big[\,A_+{U^+}^2\big|\phi_{2}^++\overline {\phi_{0}^+}\big|^2+A_-{U^-}^2\big|\phi_{2}^-+\overline {\phi_{0}^-}\big|^2\,\Big]r\,{\mathrm d}r
\nonumber\\[1mm]
&+\int_0^\infty 2BU^+U^-\Big[\,{\mathbf {Re}}\big((\phi_{2}^++\overline {\phi_{0}^+})\overline {\phi_{2}^-}\,\big)
+{\mathbf {Re}}\big((\phi_{0}^++\overline {\phi_{2}^+})\overline {\phi_{0}^-}\,\big)\,\Big]r\,{\mathrm d}r
\nonumber\\[1mm]
&+\int_0^\infty \Big[\,A_+({U^+}^2-{t^+}^2)+B({U^-}^2-{t^-}^2)\,\Big](|\phi_{2}^+|^2+|\phi_{0}^+|^2)r\,{\mathrm d}r
\nonumber\\[1mm]
&+\int_0^\infty \Big[\,A_-({U^-}^2-{t^-}^2)+B({U^+}^2-{t^+}^2)\,\Big](|\phi_{2}^-|^2+|\phi_{0}^-|^2)r\,{\mathrm d}r.
\end{align}
For any $\phi\in \mathcal S$, we can write it as  $\phi=\phi^\mathbb{R}+i\phi^\mathbb{I}$, and then get
\begin{align}\label{1.4}
\mathbb B_1 (\phi_{2},\phi_{0})
=\mathbb B_1 (\phi_{2}^\mathbb{R},\phi_{0}^\mathbb{R})+\mathbb B_1 (i\phi_{2}^\mathbb{I},i\phi_{0}^\mathbb{I})
={\mathbb D}(-\phi_{2}^\mathbb{R},\phi_{0}^\mathbb{R})+ {\mathbb D} (\phi_{2}^\mathbb{I},\phi_{0}^\mathbb{I}),
\end{align}
 with the quadratic form ${\mathbb D}$  defined as
 $$
 {\mathbb D} (\varpi, \hbar)=\mathbb B_1 (i\varpi, i\hbar),~\text{for any}~ (\varpi, \hbar)\in {\tilde{\mathcal S}}\times {\tilde{\mathcal S}}.
 $$
Here we  denote by ${\tilde{\mathcal S}}$ the subspace of real-valued functions in $\mathcal S$.
The relation in (\ref{1.4}) implies that we shall concern the quadratic form $\mathbb D$ on ${\tilde{\mathcal S}}\times{\tilde{\mathcal S}}$.

\medskip
We pause here to give an observation of the null space of $\mathbb D$.
By recalling
\begin{align}&\frac{\partial w^\pm}{\partial{x_1}}=\frac{1}{2}\Big({U^\pm}'-\frac{U^\pm}{r}\Big)e^{2i\theta}+\frac{1}{2}\Big({U^\pm}'+\frac{U^\pm}{r}\Big),
\label{1.6}
\vspace{3mm}\\
 &\frac{\partial w^\pm}{\partial{x_2}}=\frac{-i}{2}\Big({U^\pm}'-\frac{U^\pm}{r}\Big)e^{2i\theta}+\frac{i}{2}\Big({U^\pm}'+\frac{U^\pm}{r}\Big),
\label{1.60}
\end{align}
we define
\begin{align}\label{1.7}
\Phi_2=(\Phi_2^+,\Phi_2^-),~~~\Phi_0=(\Phi_0^+,\Phi_0^-),
\end{align}
with
\begin{align}\label{1.5}
\Phi_2^\pm= \frac{1}{2}\Big({-U^\pm}'+\frac{U^\pm}{r}\Big),\qquad \Phi_0^\pm=\frac{1}{2}\Big({U^\pm}'+\frac{U^\pm}{r}\Big).
\end{align}
Using the result (\ref{1.4}) together with (\ref{1.6})-(\ref{1.5}), we can obtain
\begin{align*}
&0=\frac{1}{\pi}\boldsymbol{B}\Big(\frac{\partial w}{\partial{x_1}}, \frac{\partial w}{\partial{x_1}}\Big)
= \mathbb B_1 ( -\Phi_2, \Phi_0)
= {\mathbb D} ( \Phi_2, \Phi_0),
\\[1mm]
&0=\frac{1}{\pi}\boldsymbol{B}\Big(\frac{\partial w}{\partial{x_2}}, \frac{\partial w}{\partial{x_2}}\Big)
= \mathbb B_1 ( i\Phi_2, i\Phi_0)
= {\mathbb D} ( \Phi_2, \Phi_0).
\end{align*}
Thus we have
\begin{align}\label{1.11}
{\mathbb D} (\Phi_2,\Phi_0)=0.
\end{align}

\medskip
The properties of ${\mathbb D} (\phi_2,\phi_0)$ on ${\tilde{\mathcal S}}\times {\tilde{\mathcal S}}$ will
be provided in the following proposition, which confirms the non-negativity of the quadratic form $\mathbb B_1$ on ${\mathcal S}\times {\mathcal S}$
due to the relation \eqref{1.4}.

\begin{proposition}\label{proposition2.3}
We have
$$
{\mathbb D} (\phi_2,\phi_0)\geq 0,
\quad\forall\, (\phi_2,\phi_0)\in {\tilde{\mathcal S}}\times {\tilde{\mathcal S}}
$$
where the equality holds if and only if $(\phi_2, \phi_0)=C(\Phi_2,\Phi_0)$ with $(\Phi_2, \Phi_0)$ in (\ref{1.7})-(\ref{1.5}).
\end{proposition}

{\textit{Proof of Proposition \ref{proposition2.3}}}: For any $\phi_2, \phi_{0} \in  {\tilde{\mathcal S}}\times {\tilde{\mathcal S}}$,
we compute ${\mathbb D} (\phi_2,\phi_0)$ in the following
\begin{align}\label{1.9}
&{\mathbb D} (\phi_2,\phi_0)=\mathbb B_1 (i\phi_{2},i\phi_{0})=\mathbb B_1 (-\phi_{2},\phi_{0})
\nonumber\\[1mm]
=
&\int_0^\infty \Big[\,|(\phi_{2}^+)'|^2+|(\phi_{0}^+)'|^2+\frac{4}{r^2}|\phi_{2}^+|^2\,\Big]r\,{\mathrm d}r
\,+\,\int_0^\infty \Big[\,|(\phi_{2}^-)'|^2+|(\phi_{0}^-)'|^2+\frac{4}{r^2}|\phi_{2}^-|^2\,\Big]r\,{\mathrm d}r
\nonumber\\[1mm]
&\,+\,\int_0^\infty \Big[\,A_+{U^+}^2|\phi_{0}^+-\phi_{2}^+|^2+A_-{U^-}^2|\phi_{0}^--\phi_{2}^-|^2\,\Big]r\,{\mathrm d}r
\nonumber\\[1mm]
&\,+\,\int_0^\infty 2BU^+U^-\Big[\,{\mathbf {Re}}\big((\phi_{0}^+-\phi_{2}^+){\color{blue}(-\phi_{2}^-)}\big)
+{\mathbf {Re}}\big((\phi_{0}^+-\phi_{2}^+){\color{blue}\phi_{0}^-}\big)\,\Big]r\,{\mathrm d}r\nonumber
\\[1mm]
&\,+\,\int_0^\infty \Big[\,A_+\big({U^+}^2-{t^+}^2\big)+B\big({U^-}^2-{t^-}^2\big)\,\Big](|\phi_{2}^+|^2+|\phi_{0}^+|^2)r\,{\mathrm d}r
\nonumber\\[1mm]
&\,+\,\int_0^\infty \Big[\,A_-\big({U^-}^2-{t^-}^2\big)+B\big({U^+}^2-{t^+}^2\big)\,\Big](|\phi_{2}^-|^2+|\phi_{0}^-|^2)r\,{\mathrm d}r.
\end{align}

\medskip
Furthermore, we denote
\begin{align}\label{1.8}
C^+=\phi_{0}^++\phi_{2}^+,
\quad
C^-=\phi_{0}^-+\phi_{2}^-,
\quad
D^+=\phi_{0}^+-\phi_{2}^+,
\quad
D^-=\phi_{0}^--\phi_{2}^-.
\end{align}
By substituting (\ref{1.8}) in (\ref{1.9}), we can obtain
\begin{align}\label{1.10}
{\mathbb D} (\phi_2,\phi_0)
=&\int_0^\infty\frac{1}{2}\Big[\,(C^+)'^2+(D^+)'^2+\frac{2}{r^2}(C^+-D^+)^2\,\Big]r\,{\mathrm d}r
\nonumber\\[1mm]
&\,+\,\int_0^\infty\frac{1}{2}\Big[\,(C^-)'^2+(D^-)'^2+\frac{2}{r^2}(C^--D^-)^2\,\Big]r\,{\mathrm d}r
\nonumber\\[1mm]
&\,+\,\int_0^\infty\big(A_+{U^+}^2{D^+}^2+A_-{U^-}^2{D^-}^2\big)r\,{\mathrm d}r
\nonumber\\[1mm]
&\,+\,\int_0^\infty\frac{1}{2}\Big[\,A_+\big({U^+}^2-{t^+}^2\big)+B\big({U^-}^2-{t^-}^2\big)\,\Big]({C^+}^2+{D^+}^2)r\,{\mathrm d}r
\nonumber\\[1mm]
&\,+\,\int_0^\infty\frac{1}{2}\Big[\,A_-\big({U^-}^2-{t^-}^2\big)+B\big({U^+}^2-{t^+}^2\big)\,\Big]({C^-}^2+{D^-}^2)r\,{\mathrm d}r
\nonumber\\[1mm]
&\,+\,\int_0^\infty 2B U^+U^-D^+D^-r\mathrm \,{\mathrm d}r
\nonumber\\[1mm]
\equiv&\,F(C^+,C^-,D^+,D^-).
\end{align}
On the other hand, we can write the functional  $F({{\mathbf u}^+},{{\mathbf v}^+},{{\mathbf u}^-},{{\mathbf v}^-})$ in the form
\begin{align*}
F({{\mathbf u}^+},{{\mathbf v}^+},{{\mathbf u}^-},{{\mathbf v}^-})=&
\int_0^\infty\Big[\,\alpha_1({{\mathbf u}^+}')^2+\alpha_2({{\mathbf u}^-}')^2+\beta_1({{\mathbf v}^+}')^2+\beta_2({{\mathbf v}^-}')^2\,\Big]
\\[1mm]
&\,+\,\int_0^\infty \Big[\,a_1({{\mathbf u}^+})^2+a_2({{\mathbf u}^-})^2+d_1({{\mathbf v}^+})^2+d_2({{\mathbf v}^-})^2\,\Big]
\\[1mm]
&\,+\,\int_0^\infty \Big[\,2b_1{\mathbf u}^+{\mathbf v}^++2b_2{\mathbf u}^-{\mathbf v}^-+2b_3{\mathbf v}^-{\mathbf v}^+\,\Big],
\end{align*}
where
$$ a_1=\frac{1}{r}+\frac{1}{2}\Big[\,A_+\big({U^+}^2-{t^+}^2\big)+B\big({U^-}^2-{t^-}^2\big)\,\Big]r,$$
$$ a_2=\frac{1}{r}+\frac{1}{2}\Big[\,A_-\big({U^-}^2-{t^-}^2\big)+B\big({U^+}^2-{t^+}^2\big)\,\Big]r,$$
$$ d_1=\frac{1}{r}+\frac{1}{2}\Big[\,A_+\big({U^+}^2-{t^+}^2\big)+B\big({U^-}^2-{t^-}^2\big)\,\Big]r+A_+{U^+}^2r,$$
$$ d_2=\frac{1}{r}+\frac{1}{2}\Big[\,A_-\big({U^-}^2-{t^-}^2\big)+B\big({U^+}^2-{t^+}^2\big)\,\Big]r+A_-{U^-}^2r.$$
In the above, we also have denoted
\begin{align*}
&\alpha_1=\frac{r}{2},\quad \alpha_2=\frac{r}{2},\quad \beta_1=\frac{r}{2},\quad \beta_2=\frac{r}{2},
\\[1mm]
&b_1=\frac{-1}{r},\quad b_2=\frac{-1}{r},\quad b_3=B U^+U^-r.
\end{align*}
Note that,  for all $r>0$
\begin{align}\label{negativityb1b2b3}
 \alpha_1, \alpha_2, \beta_1, \beta_2>0
\quad  \mbox{and}\quad
b_1, b_2, b_3<0,
\end{align}
under the assumption $B<0$ in (H2).

\medskip
According to \eqref{1.8}, for the functions $\Phi_0$ and $\Phi_2$ in (\ref{1.7})-(\ref{1.5}), we set
\begin{align}
\eta^+=\Phi_{0}^++\Phi_{2}^+={\frac{U^+}{r}},\qquad\eta^-=\Phi_{0}^-+\Phi_{2}^-={\frac{U^-}{r}},\label{1.42}
\\[1mm]
\zeta^+=\Phi_{0}^+-\Phi_{2}^+={U^+}',\qquad\zeta^-=\Phi_{0}^--\Phi_{2}^-={U^-}'.\label{1.43}
\end{align}
A direct computation together with the equations (\ref{1.2}) shows that  $(\eta^+, \zeta^+, \eta^-,  \zeta^-)$ satisfies the equations
\begin{align}\label{1.14}
\begin{cases}
 \big(-\alpha_1{\eta^+}'\big)'\,+\,a_1 \eta^+\,+\,b_1\zeta^+\,=\,0,
 \\[2mm]
 \big(-\alpha_2{\eta^-}'\big)'\,+\,a_1 \eta^-\,+\,b_2\zeta^-\,=\,0,
 \\[2mm]
 \big(-\beta_1{\zeta^+}'\big)'\,+\,d_1 \zeta^+\,+\,b_1\eta^+\,+\,b_3\zeta^-\,=\,0,
 \\[2mm]
 \big(-\beta_2{\zeta^-}'\big)'\,+\,d_2 \zeta^-\,+\,b_2\eta^-\,+\,b_3\zeta^+\,=\,0.
 \end{cases}
\end{align}
For any smooth functions $({{\mathbf u}^+}, {{\mathbf v}^+}, {{\mathbf v}^-}{{\mathbf u}^-}, {{\mathbf v}^-})\in C_c(0,\infty)\times C_c(0,\infty)\times C_c(0,\infty) \times C_c(0,\infty)$, the Picone's identities imply that
 \begin{align}
\big({{{\mathbf u}^+}}'\big)^2\,-\,\Big(\frac{{{{\mathbf u}^+}}^2}{\eta^+}\Big)'(\eta^+)'\,=\,\Big({{{\mathbf u}^+}}'-\frac{{{\mathbf u}^+}}{\eta^+}{\eta^+}'\Big)^2\,\geq\, 0,\label{1.13}
 \\[1mm]
\big({{{\mathbf u}^-}}'\big)^2\,-\,\Big(\frac{{{{\mathbf u}^-}}^2}{\eta^-}\Big)'(\eta^-)'\,=\,\Big({{{\mathbf u}^-}}'-\frac{{{\mathbf u}^-}}{\eta^-}{\eta^-}'\Big)^2\,\geq\, 0,
\\[1mm]
\big({{{\mathbf v}^+}}'\big)^2\,-\,\Big(\frac{{{{\mathbf v}^+}}^2}{\zeta^+}\Big)'(\zeta^+)'\,=\,\Big({{{\mathbf v}^+}}'-\frac{{{\mathbf v}^+}}{\zeta^+}{\zeta^+}'\Big)^2\,\geq\, 0,
\\[1mm]
\big({{{\mathbf v}^-}}'\big)^2\,-\,\Big(\frac{{{{\mathbf v}^-}}^2}{\zeta^-}\Big)'(\zeta^-)'\,=\,\Big({{{\mathbf v}^-}}'-\frac{{{\mathbf v}^-}}{\zeta^-}{\zeta^-}'\Big)^2\,\geq\, 0.\label{1.15}
\end{align}

\medskip
Now, we consider the properties of $F$.
Multiplying (\ref{1.13})-(\ref{1.15}) by $\alpha_1,\alpha_2,\beta_1,\beta_2$ respectively together with integrating by parts, and using (\ref{1.14}), we can get
\begin{align}\label{1.16}
0\leq &\int_0^\infty \Big[\,\alpha_1({{{\mathbf u}^+}}')^2+\alpha_2({{{\mathbf u}^-}}')^2 +\beta_1({{{\mathbf v}^+}}')^2+\beta_2({{{\mathbf v}^-}}')^2\,\Big] \,{\mathrm d}r
\nonumber\\[1mm]
&\,+\, \int_0^\infty \Big[\,a_1 {{{\mathbf u}^+}}^2+a_2 {{{\mathbf u}^-}}^2+d_1 {{{\mathbf v}^+}}^2+d_2 {{{\mathbf v}^-}}^2\,\Big]\,{\mathrm d}r
\nonumber\\[1mm]
&\,+\,\int_0^\infty \Big[\,b_1{{{\mathbf u}^+}}^2\frac{\zeta^+}{\eta^+}+b_2{{{\mathbf u}^-}}^2\frac{\zeta^-}{\eta^-}+ b_1{{{\mathbf v}^+}}^2\frac{\eta^+}{\zeta^+}+b_2{{{\mathbf v}^-}}^2\frac{\eta^-}{\zeta^-}+b_3{{{\mathbf v}^+}}^2\frac{\zeta^-}{\zeta^+}+b_3{{{\mathbf v}^-}}^2\frac{\zeta^+}{\zeta^-}\,\Big]\,{\mathrm d}r
\nonumber\\[1mm]
=&\,F({{\mathbf u}^+},{{\mathbf v}^+},{{\mathbf u}^-},{{\mathbf v}^-})
\,+\,\int_0^\infty b_1\Big[\,{{\mathbf u}^+}\Big(\frac{\zeta^+}{\eta^+}\Big)^{\frac{1}{2}}-{{\mathbf v}^+}\Big(\frac{\eta^+}{\zeta^+}\Big)^{\frac{1}{2}}\,\Big]^2\,{\mathrm d}r
\nonumber\\[1mm]
&\,+\,
\int_0^\infty b_2\Big[\,{{\mathbf u}^-}\Big(\frac{\zeta^-}{\eta^-}\Big)^{\frac{1}{2}}-{{\mathbf v}^-}\Big(\frac{\eta^-}{\zeta^-}\Big)^{\frac{1}{2}}\,\Big]^2\,{\mathrm d}r
\,+\,\int_0^\infty b_3\Big[\,{{\mathbf v}^+}\Big(\frac{\zeta^-}{\zeta^+}\Big)^{\frac{1}{2}}-{{\mathbf v}^-}\Big(\frac{\zeta^+}{\zeta^-}\Big)^{\frac{1}{2}}\,\Big]^2\,{\mathrm d}r.
\end{align}
This implies that
\begin{align}\label{1.41}
F({{\mathbf u}^+},{{\mathbf v}^+},{{\mathbf u}^-},{{\mathbf v}^-})\geq
&\,-\,\int_0^\infty b_1\Big[\,{{\mathbf u}^+}(\frac{\zeta^+}{\eta^+})^{\frac{1}{2}}-{{\mathbf v}^+}(\frac{\eta^+}{\zeta^+})^{\frac{1}{2}}\,\Big]^2
\,-\,\int_0^\infty b_2\Big[\,{{\mathbf u}^-}(\frac{\zeta^-}{\eta^-})^{\frac{1}{2}}-{{\mathbf v}^-}(\frac{\eta^-}{\zeta^-})^{\frac{1}{2}}\,\Big]^2
\nonumber\\[1mm]
&\,-\,\int_0^\infty b_3\Big[\,{{\mathbf v}^+}(\frac{\zeta^-}{\zeta^+})^{\frac{1}{2}}-{{\mathbf v}^-}(\frac{\zeta^+}{\zeta^-})^{\frac{1}{2}}\,\Big]^2
\nonumber\\[1mm]
\geq&\, 0,
\end{align}
due to the facts in (\ref{negativityb1b2b3}).
Inequalities in (\ref{1.41}) also hold
for $({{\mathbf u}^+}, {{\mathbf v}^+}, {{\mathbf u}^-}, {{\mathbf v}^-})\in {\tilde{\mathcal S}}\times {\tilde{\mathcal S}}\times{\tilde{\mathcal S}}\times {\tilde{\mathcal S}}$ due to a density argument. The inequalities in (\ref{1.41}) imply  $F({{\mathbf u}^+},{{\mathbf v}^+},{{\mathbf u}^-},{{\mathbf v}^-})\geq 0$.
Moreover, if the equality $F({{\mathbf u}^+},{{\mathbf v}^+},{{\mathbf u}^-},{{\mathbf v}^-})=0$ holds, then there exists a function $\vartheta$ such that
$$
({{\mathbf u}^+},{{\mathbf v}^+},{{\mathbf u}^-},{{\mathbf v}^-})=(\vartheta\eta^+,\vartheta\zeta^+,\vartheta\eta^-,\vartheta\zeta^-).
 $$
Applying (\ref{1.13})-(\ref{1.15}), we know that the function $\vartheta$  must be a constant.  Therefore, we conclude that  $F({{\mathbf u}^+},{{\mathbf v}^+},{{\mathbf u}^-},{{\mathbf v}^-})\geq 0$ and the equality holds if and only if
$$
({{\mathbf u}^+},{{\mathbf v}^+},{{\mathbf u}^-},{{\mathbf v}^-})=C(\eta^+,\zeta^+,\eta^-,\zeta^-).
$$

\medskip
Combining with (\ref{1.8}), (\ref{1.10}), (\ref{1.42}), (\ref{1.43}) and the results above, we have
\begin{align*}
{\mathbb D} (\phi_2,\phi_0)\geq 0,
\end{align*}
and if $ {\mathbb D} (\phi_2,\phi_0)= 0$, then
$$(\phi_2,\phi_0)=C(\Phi_2,\Phi_0).$$
This finishes  the proof of the Proposition \ref{proposition2.3}\vspace{3mm}.
\qed

\medskip
\begin{remark}\label{02}
In the proof of Proposition \ref{proposition2.3}, in order to get the result
$$
F({{\mathbf u}^+},{{\mathbf v}^+},{{\mathbf u}^-},{{\mathbf v}^-})\geq 0,
$$
we use the facts in (\ref{negativityb1b2b3}).
That is why we need the technical assumption $B<0$ in Theorem \ref{theorem1.1}.
For the case $B=0$, the coupled Ginzburg-Landau system (\ref{1.1}) is reduced to the single-component  Ginzburg-Landau equation in $\mathbb R^2$, and the non-degeneracy results have been proved in \cite{DFK}. When $B>0$ the non-degeneracy type result
for local energy minimizer of the system (\ref{1.1}) is not clear.
\qed
\end{remark}

\medskip
Now we give the proof of  Theorem \ref{theorem1.1}.\\
{\textit{Proof of Theorem \ref{theorem1.1}}}: For any $\phi\in \mathcal H$, we can decompose $\boldsymbol{B}(\phi,\phi)$ as in (\ref{1.45}).
By combining Propositions \ref{proposition2.1}-\ref{proposition2.3} and (\ref{1.4}),
we can know that $\boldsymbol{B}(\phi,\phi)\geq 0$.
In particular, by the equalities in Propositions \ref{proposition2.1}-\ref{proposition2.3},
we have that
$$
\boldsymbol{B}(\phi,\phi)= 0\quad \text{iff}\quad \phi= \phi_0+\phi_1e^{i\theta}+\phi_2e^{2i\theta},
$$
 where
\begin{align*}
(\phi_2^{\mathbb I}, \phi_0^{\mathbb I})=c_1(\Phi_2,\, \Phi_0),
\qquad
(-\phi_2^{\mathbb R}, \phi_0^{\mathbb R})=c_2(\Phi_2,\, \Phi_0),
\qquad
\phi_1= i (c_3\, U^+, c_4\,U^-),
\end{align*}
for some constants $c_1, c_2, c_3, c_4\in \mathbb{R}$.
It is easy to verify that
\begin{align*}
 {\mathcal L}(\phi)=0  \quad\text{iff}\quad
\phi=c_1\frac{\partial w}{\partial x_1}+c_2\frac{\partial w}{\partial x_2}+i(c_3 \,w^+,c_4\,w^-)~~~\text{for any } \phi \in \mathcal H.
\end{align*}
Since $i(c_3\,w^+,c_4\,w^-) $ is not contained in $\mathcal H$, we have
\begin{align*}
\phi=c_1\frac{\partial w}{\partial x_1}+c_2\frac{\partial w}{\partial x_2}.
\end{align*}
This completes  the proof of Theorem \ref{theorem1.1}.
\qed

\medskip
\section{Application of Non-Degeneracy: Fredholm Alternative.}\label{section3}

In the section,   we will find a solution $\psi=(\psi^+, \psi^-): {\mathbb R}^2\rightarrow {\mathbb C}^2$ of the linear equation
\begin{eqnarray}\label{1.23}
{\mathcal L}(\psi)= h,
%\left\{\begin{array}{rcl}  {\mathcal L}_+(\psi)=h^+,
%&&\\
%    \space{ }\\
%  {\mathcal L}_-(\psi)=h^-.\end{array}\right,
\end{eqnarray}
for a given $h=(h^+, h^-): {\mathbb R}^2\rightarrow {\mathbb C}^2$.
The idea is to find a minimizer of energy functional corresponding to the equation
(\ref{1.23})
\begin{align}\label{functionalJ}
J(\psi)=\frac{1}{2}\boldsymbol{B}(\psi,\psi)-\langle \psi,h\rangle,
\end{align}
where $\boldsymbol{B}(\cdot, \cdot)$ and $\langle \cdot,\cdot\rangle$ are defined in (\ref{0.2}) and (\ref{0.3}).

\medskip
\subsection{\textbf{ Preliminary work for the proof of Theorem \ref{theorem1.2} }}

We give some preliminaries for the proof of Theorem \ref{theorem1.2}.
We decompose $\psi=({\psi}^+, {\psi}^-)$ into the form
\begin{align}\label{decompositionphi}
\psi=\psi_0+\sum_{j=1}^\infty{\psi_j^1}+\sum_{j=1}^\infty{\psi_j^2},
\end{align}
where
\begin{align}\label{phi0}
\psi_0=({\psi_0}^+,{\psi_0}^-),\quad \psi_j^1=({\psi_j^1}^+,{\psi_j^1}^-),\quad \psi_j^2=({\psi_j^2}^+,{\psi_j^2}^-),
\end{align}
 and
 \begin{equation}\label{decompositionh0}
  \begin{split}
 &{\psi_0}^{\pm}=e^{i\theta}\Big[\,{\psi_{01}}^{\pm}\,+\,i{\psi_{02}}^{\pm}\,\Big],
 \\[1mm]
  &{\psi_j^1}^{\pm}=e^{i\theta}\Big[\,{\psi_{j1}^1}^{\pm}\sin{j\theta}\,+\,i{\psi_{j2}^1}^{\pm}\cos{j\theta}\,\Big],
  \\[1mm]
  &{\psi_j^2}^{\pm}=e^{i\theta}\Big[\,{\psi_{j1}^2}^{\pm}\cos{j\theta}\,+\,i{\psi_{j2}^2}^{\pm}\sin{j\theta}\,\Big].
  \end{split}
 \end{equation}
  The decompositions in (\ref{decompositionphi})-(\ref{decompositionh0})  are naturally associated to the functions
 $ \frac{\partial w}{\partial x_1}, \frac{\partial w}{\partial x_2}$ in the kernel  of the linearized operator ${\mathcal L}$.
 In fact, we can write
 \begin{align}
 &\frac{\partial w}{\partial x_1}=e^{i\theta} \Big[U'\cos \theta-i \frac{U}{r}\sin \theta\,\Big],\label{iw}
 \\[1mm]
 & \frac{\partial w}{\partial x_2}=e^{i\theta} \Big[U'\sin \theta+i \frac{U}{r}\cos \theta\,\Big].\label{iw2}
   \end{align}
\medskip
Similarly, we  decompose $h=({h}^+,{h}^-)$
\begin{align}\label{decompositionh}
h=h_0+\sum_{j=1}^\infty{h_j^1}+\sum_{j=1}^\infty{h_j^2},
\end{align}
where
$$
h_0=({h_0}^+,  {h_0}^-),\quad h_j^1=({h_j^1}^+, {h_j^1}^-),\quad h_j^2=({h_j^2}^+, {h_j^2}^-),
$$
and
\begin{align*}
 &{h_0}^{\pm}=e^{i\theta}\Big[\,{h_{01}}^{\pm}\,+\,i{h_{02}}^{\pm}\,\Big],
 \\[1mm]
  &{h_j^1}^{\pm}=e^{i\theta}\Big[\,{h_{j1}^1}^{\pm}\sin{j\theta}\,+\,i{h_{j2}^1}^{\pm}\cos{j\theta}\,\Big],
  \\[1mm]
  &{h_j^2}^{\pm}=e^{i\theta}\Big[\,{h_{j1}^2}^{\pm}\cos{j\theta}\,+\,i{h_{j2}^2}^{\pm}\sin{j\theta}\,\Big].
 \end{align*}
Then finding a solution to (\ref{1.23}) is equivalent to solving the equations
\begin{align}
&{\mathcal L}(\psi_0)=h_0,\label{1.47}
\\[1mm]
&{\mathcal L}(\psi_j^\ell)=h_j^\ell,~~~j\in \mathbb{N}^+,\, \ell=1,2.\label{1.31}
\end{align}
We will solve equations (\ref{1.47}) and (\ref{1.31}) separately and provide the
proof for Theorem \ref{theorem1.2} in three steps, see Section \ref{section3.2}.
For that purpose, it is natural to consider the critical point of the functionals

\begin{equation}\label{functionalJ2}
J_0(\psi_0)=\frac{1}{2}\boldsymbol{B} (\psi_0, \psi_0)-\langle \psi_0,h_0\rangle,
\end{equation}
\begin{equation}\label{functionalJ3}
J_j^\ell(\psi_j^\ell)=\frac{1}{2}\boldsymbol{B}(\psi_j^\ell, \psi_j^\ell)-\langle \psi_j^\ell,h_j^\ell\rangle,~~~j\in \mathbb{N}^+,\, \ell=1,2.
\end{equation}

\medskip
\subsection{\textbf{Proof of Theorem \ref{theorem1.2} }}\label{section3.2}

As we have said in previous section, in order to prove Theorem \ref{theorem1.2}, we will solve the equations in (\ref{1.47}) and (\ref{1.31}).
\\

\noindent{\textbf{ The mode $j=1$.}}
We now  begin by solving the problem ${\mathcal L}(\psi_1^1)=h_1^1$.
Let $\mathcal{H_*}$ be the space of functions
$$
{\tilde\psi}=({\tilde\psi}^+, {\tilde\psi}^-)=\Big(({\tilde\psi}_1^+, {\tilde\psi}_2^+), ({\tilde\psi}_1^-, {\tilde\psi}_2^-)\Big)
$$
such that the functions $\psi=({\psi}^{+},{\psi}^{-})$ belong to $\mathcal H$ with the forms
\begin{align}\label{phi1}
{\psi}^{\pm}=e^{i\theta}\big[{\tilde\psi_1}^{\pm}\sin{\theta}+i{\tilde\psi_2}^{\pm}\cos{\theta}\big].
\end{align}
According to the norm for $\psi$ given in \eqref{mathcalH}, the norm for ${\tilde\psi}$ in the space $\mathcal{H_*}$  is given by
\begin{align*}
\|\tilde\psi\|^2_{\mathcal{H_*}}=&\int_0^\infty \Big[\,\big|(\tilde\psi^+)'\big|^2+\frac{1}{r^2}\big(\tilde\psi_1^+-\tilde\psi_2^+\big)^2\,\Big]r\,{\mathrm d}r
+\int_0^\infty \Big[\,\big|(\tilde\psi^-)'\big|^2+\frac{1}{r^2}\big(\tilde\psi_1^--\tilde\psi_2^-\big)^2\,\Big]r\,{\mathrm d}r
\nonumber\\[1mm]
&+\int_0^\infty \Big[\,A_+\big({t^+}^2-{U^+}^2\big)-B\big({t^-}^2-{U^-}^2\big)\,\Big]|\tilde\psi^+|^2r\,{\mathrm d}r
\nonumber\\[1mm]
&+\int_0^\infty \Big[\,A_-\big({t^-}^2-{U^-}^2\big)-B\big({t^+}^2-{U^+}^2\big)\,\Big]|\tilde\psi^-|^2r\,{\mathrm d}r.
\end{align*}
%Next, we introduce  the new   quadratic form  $\mathscr{B}(\tilde\psi,\tilde\psi)$ for the function $ \tilde\psi \in  \mathcal{H_*}$.

\medskip
We  introduce  new   quadratic form  $\mathscr{B}(\tilde\psi,\tilde\psi)$ for the function $ \tilde\psi \in  \mathcal{H_*}$,  which is defined as
\begin{align}\label{mathcal{B}_1^1}
\mathscr{B}(\tilde\psi,\tilde\psi) \,=\, \pi\, \boldsymbol{B}(\psi,\psi),
\end{align}
where $\psi=({\psi}^{+},{\psi}^{-})$  is the function  defined in \eqref{phi1}.
By the relations in \eqref{mathcal{B}_1^1}, we can write quadratic form  $\mathscr{B}(\tilde\psi,\tilde\psi)$ as follows
\begin{align}\label{mathscr b phi}
\mathscr{B}(\tilde\psi,\tilde\psi)
=&\int_0^\infty \Big[\,\big|(\tilde\psi^+)'\big|^2+\frac{2}{r^2}\big(\tilde\psi_1^+-\tilde\psi_2^+\big)^2\,\Big]r\,{\mathrm d}r
+\int_0^\infty \Big[\,\big|(\tilde\psi^-)'\big|^2+\frac{2}{r^2}\big(\tilde\psi_1^--\tilde\psi_2^-\big)^2\,\Big]r\,{\mathrm d}r
\nonumber\\[1mm]
&+\int_0^\infty \Big[\,A_+\big({U^+}^2-{t^+}^2\big)+B\big({U^-}^2-{t^-}^2\big)\,\Big]|\tilde\psi^+|^2r\,{\mathrm d}r
\nonumber\\[1mm]
&+\int_0^\infty \Big[\,A_-\big({U^-}^2-{t^-}^2\big)+B\big({U^+}^2-{t^+}^2\big)\,\Big]|\tilde\psi^-|^2r\,{\mathrm d}r
\nonumber\\[1mm]
&+\int_0^\infty\Big[\,2A_+{U^+}|\tilde\psi_1^+|^2+2A_-{U^-}^2|\tilde\psi_1^-|^2+4BU^+U^-\tilde\psi_1^+\tilde\psi_1^-\,\Big]r\,{\mathrm d}r.
\end{align}

Define
$$
\tilde h=(\tilde h^+, \tilde h^-)=\Big(\big({h_{11}^1}^+, {h_{12}^1}^+\big), \big({h_{11}^1}^-, {h_{12}^1}^-\big)\Big).
$$
Then, for $\ell=1,j=1$, solving problem (\ref{1.31})   in  $\mathcal H$ is corresponding exactly to finding a critical point of the functional $J_1(\tilde\psi)$, which is given by
\begin{align*}
J_1(\tilde\psi)=\frac{1}{2}\mathscr{B}(\tilde\psi,\tilde\psi)
-\int_0^\infty\big({h_{11}^1}^+{\tilde\psi}_1^++{h_{11}^1}^-{\tilde\psi}_1^-+{h_{12}^1}^+{\tilde\psi}_2^+
+{h_{12}^1}^-{\tilde\psi}_2^-\big)r\,{\mathrm d}r.
\end{align*}

\medskip
Now denote
$$
Z_0=(Z_0^+, Z_0^-)
\qquad\mbox{with}\quad
Z_0^+=\Big({U^+}',\frac{U^+}{r}\Big),
\quad
Z_0^-=\Big({U^-}',\frac{U^-}{r}\Big).
 $$
Combining  Theorem \ref{theorem1.1} with  \eqref{iw2} and \eqref{mathcal{B}_1^1}, we can know easily that
$\mathscr{B}(\tilde\psi,\tilde\psi)=0$ if and only if $\tilde\psi=CZ_0$. The assumption (\ref{1.35}) implies that
\begin{align}\label{a.1}
\int_0^\infty\Big({h_{11}^1}^+{U^+}'+{h_{11}^1}^-{U^-}'+{h_{12}^1}^+\frac{U^+}{r}
+{h_{12}^1}^-\frac{U^-}{r}\Big)r\,{\mathrm d}r=0.
\end{align}

\medskip
We define a weighted inner product $\langle \cdot,\cdot\rangle_*$ for the space $\mathcal{H_*}$
\begin{align*}
\langle u,v\rangle_*=\langle u^+,\, v^+\rangle_{*\mathcal R^+}\,+\,\langle {u^-},\, v^-\rangle_{*\mathcal R^-},
\end{align*}
for any
$$
u=(u^+,u^-)=\Big((u^+_1,u^+_2),(u^-_1,u^-_2)\Big)\in \mathcal{H_*},
\qquad
v=(v^+,v^-)=\Big((v^+_1,v^+_2),(v^-_1,v^-_2)\Big)\in \mathcal{H_*},
$$
where
\begin{align*}
\langle u^+,v^+\rangle_{*\mathcal R^+}=
\int_0^\infty \Big[\,A_+\big({t^+}^2-{U^+}^2\big)-B\big({t^-}^2-{U^-}^2\big)\,\Big](u^+_1v^+_1+u^+_2v^+_2)r\,{\mathrm d}r,
\end{align*}
\begin{align*}
\langle u^-,v^-\rangle_{*\mathcal R^-}=
\int_0^\infty \Big[\,A_-\big({t^-}^2-{U^-}^2\big)-B\big({t^+}^2-{U^+}^2\big)\,\Big](u^-_1v^-_1+u^-_2v^-_2)r\,{\mathrm d}r.
\end{align*}
In order to solve the equation (\ref{1.31}) for $j=1,\ell=1$, we need the following lemma.

\begin{lemma}\label{1.49}
There exists a constant $C>0$ such that for any
$$
\tilde\psi=(\tilde\psi^+, \tilde\psi^-)=\Big(\big({\tilde\psi}_1^+, {\tilde\psi}_2^+\big), \big({\tilde\psi}_1^-, {\tilde\psi}_2^-\big)\Big)\in \mathcal{H_*},
$$
satisfying
\begin{align}
\langle \tilde\psi, Z_0\rangle_*=\langle {\tilde\psi}^+,Z_0^+\rangle_{*\mathcal R^+}+\langle {\tilde\psi}^-,Z_0^-\rangle_{*\mathcal R^-}=0,
\end{align}
%with $$Z_0^+=\Big(-{U^+}',\frac{U^+}{r}\Big),~~~Z_0^-=\Big(-{U^-}',\frac{U^-}{r}\Big),$$
we have
\begin{align}
C\|\tilde\psi\|^2_{\mathcal H_*}\leq \mathscr{B}(\tilde\psi,\tilde\psi).\end{align}
\end{lemma}

{\textit{Proof of Lemma \ref{1.49}  }}:
 Recall  the expression  of the  quadratic form $\mathscr{B}(\tilde\psi,\tilde\psi)$ in \eqref{mathscr b phi}.
% \begin{align}
%\mathscr{B}(\tilde\psi,\tilde\psi)
%=&\int_0^\infty \Big[\,\big|(\tilde\psi^+)'\big|^2+\frac{2}{r^2}\big(\tilde\psi_1^+-\tilde\psi_2^+\big)^2\,\Big]r\,{\mathrm d}r
%+\int_0^\infty \Big[\,\big|(\tilde\psi^-)'\big|^2+\frac{2}{r^2}\big(\tilde\psi_1^--\tilde\psi_2^-\big)^2\,\Big]r\,{\mathrm d}r
%\nonumber\\[1mm]
%&+\int_0^\infty \Big[\,A_+\big({U^+}^2-{t^+}^2\big)+B\big({U^-}^2-{t^-}^2\big)\,\Big]|\tilde\psi^+|^2r\,{\mathrm d}r
%\nonumber\\[1mm]
%&+\int_0^\infty \Big[\,A_-\big({U^-}^2-{t^-}^2\big)+B\big({U^+}^2-{t^+}^2\big)\,\Big]|\tilde\psi^-|^2r\,{\mathrm d}r
%\nonumber\\[1mm]
%&+\int_0^\infty\Big[\,2A_+{U^+}|\tilde\psi_1^+|^2+2A_-{U^-}^2|\tilde\psi_1^-|^2+4BU^+U^-\tilde\psi_1^+\tilde\psi_1^-\,\Big]r\,{\mathrm d}r.
%\nonumber
%\end{align}
Note that,
$$
{t^\pm}^2-{U^\pm}^2
\sim\frac{c_\pm t_\pm}{r^2}, \quad ~\text{as}~ r\rightarrow\infty,
\quad\mbox{with}\quad
c_\pm=\frac{A_\mp-B}{(A_+A_--B^2)t^\pm}.
$$
For any given $\delta>0$ small, there exists an $R>0$ large such that, for $r>R$
\begin{align}\label{1.25}
I^+(\tilde\psi)
\,=\,&\,\frac{2-\delta}{r^2}\Big(\tilde\psi_1^+-\tilde\psi_2^+\Big)^2
\,+\,2A_+{U^+}^2|\tilde\psi_1^+|^2
\,+\,2BU^+U^-\tilde\psi_1^+\tilde\psi_1^-
\nonumber\\[1mm]
&\,-\,\Big[\,A_+\big({t^+}^2-{U^+}^2\big)+B\big({t^-}^2-{U^-}^2\big)\,\Big]|\tilde\psi^+|^2
\nonumber\\[1mm]
&\,-\,\frac{\delta}{2\tau^+} \Big[\,A_+\big({t^+}^2-{U^+}^2\big)-B\big({t^-}^2-{U^-}^2\big)\,\Big]|\tilde\psi^+|^2
\nonumber\\[1mm]
\,\geq\, &\frac{1-2\delta}{r^2}|\tilde\psi^+|^2
\,+\,2A_+{U^+}^2|\tilde\psi_1^+|^2
\,+\,2BU^+U^-\tilde\psi_1^+\tilde\psi_1^-
\,-\, 2\frac{2-\delta}{r^2}{\tilde\psi}_1^+{\tilde\psi}_2^+,
\end{align}
and also
\begin{align}\label{1.26}
I^-(\tilde\psi)
\,=\,&\,\frac{2-\delta}{r^2}\Big(\tilde\psi_1^--\tilde\psi_2^-\Big)^2
\,+\,2A_-\,{U^-}^2|\tilde\psi_1^-|^2
\,+\,2B\,U^+U^-\tilde\psi_1^+\tilde\psi_1^-
\nonumber\\[1mm]
&
-\Big[\,A_-\big({t^-}^2-{U^-}^2\big)+B\big({t^+}^2-{U^+}^2\big)\,\Big]|\tilde\psi^-|^2
\nonumber\\[1mm]
&-\frac{\delta}{2\tau^-}  \Big[\,A_-\big({t^-}^2-{U^-}^2\big)-B\big({t^+}^2-{U^+}^2\big)\,\Big]|\tilde\psi^-|^2
\nonumber\\[1mm]
\,\geq\,&\, \frac{1-2\delta}{r^2}|\tilde\psi^-|^2
\,+\,2A_-\,{U^-}^2|\tilde\psi_1^-|^2
\,+\,2B\,U^+U^-\tilde\psi_1^+\tilde\psi_1^-
\,-\, 2\frac{2-\delta}{r^2}{\tilde\psi}_1^-{\tilde\psi}_2^-,
\end{align}
where
$$
\tau^\pm=\max\big\{\, A_\pm c_\pm t^\pm, |B|c_\mp t^\mp\, \big\}.
$$
Since $ A_+A_--B^2>0$,
then there exist $\Lambda>0$ such that, for $r>R$
\begin{align}
\begin{aligned}\label{1.27}
2A_+{U^+}^2|\tilde\psi_1^+|^2
\,+\,4BU^+U^-\tilde\psi_1^+\tilde\psi_1^-
\,+\,2A_-{U^-}^2|\tilde\psi_1^-|^2
\,\geq\, \Lambda \big(|\tilde\psi_1^+|^2+|\tilde\psi_1^-|^2\big).
\end{aligned}
\end{align}
%Denote
%\begin{align}
%I^+(\tilde\psi)=&\,\frac{2-\delta}{r^2}\Big(\tilde\psi_1^+-\tilde\psi_2^+\Big)^2
%\,+\,(2A_+-\delta){U^+}^2|\tilde\psi_1^+|^2
%\,+\,(2B-\delta)U^+U^-\tilde\psi_1^+\tilde\psi_1^-
%\nonumber\\[1mm]
%&\,-\,\Big[\,A_+\big({t^+}^2-{U^+}^2\big)+B\big({t^-}^2-{U^-}^2\big)\,\Big]|\tilde\psi^+|^2
%\nonumber\\[1mm]
%&\,-\,\frac{\delta}{2\tau^+} \Big[\,A_+\big({t^+}^2-{U^+}^2\big)-B\big({t^-}^2-{U^-}^2\big)\,\Big]|\tilde\psi^+|^2
%\label{1.33},
%\end{align}
%and
%\begin{align}
%I^-(\tilde\psi)=&\,\frac{2-\delta}{r^2}\Big(\tilde\psi_1^--\tilde\psi_2^-\Big)^2
%\,+\,(2A_--\delta){U^-}^2|\tilde\psi_1^-|^2
%\,+\,(2B-\delta)U^+U^-\tilde\psi_1^+\tilde\psi_1^-
%\nonumber\\[1mm]
%&\,-\,\Big[\,A_-\big({t^-}^2-{U^-}^2\big)+B\big({t^+}^2-{U^+}^2\big)\,\Big]|\tilde\psi^-|^2
%\nonumber\\[1mm]
%&\,-\,\frac{\delta}{2\tau^-} \Big[\,A_-\big({t^-}^2-{U^-}^2\big)-B\big({t^+}^2-{U^+}^2\big)\,\Big]|\tilde\psi^-|^2.
%\label{1.34}
%\end{align}
This implies that, for given $\delta, R$, when $r>R$, we have
\begin{align}
I^+(\tilde\psi)+I^-(\tilde\psi)
& \geq\frac{1-2\delta}{r^2}\Big(|\tilde\psi^+|^2+|\tilde\psi^-|^2\Big)
\,-\, 2\frac{2-\delta}{r^2}\Big(|\tilde\psi_1^+|\,|\tilde \psi_2^+|+|\tilde\psi_1^-|\,|\tilde \psi_2^-|\Big)
\nonumber
\\[1mm]
&\quad\,+\, \Lambda \Big(|\tilde\psi_1^+|^2+|\tilde\psi_1^-|^2\Big)
\nonumber
\\[1mm]
&\geq \frac{\varsigma_1}{r^2}|\tilde\psi^+|^2
\,+\,\frac{\varsigma_2}{r^2}|\tilde\psi^-|^2\geq 0,
 \end{align}
 with $\varsigma_1, \varsigma_2$ small and independent of $\tilde \psi$.
 Then there exist two positive constants $C_1, C_2$ such that
 \begin{align}\label{1.32}
 \mathscr{B}(\tilde\psi,\tilde\psi)\geq C_1\|\tilde\psi\|_{\mathcal H_*}^2-C_2\int_0^R\Big(|\tilde\psi^+|^2+|\tilde\psi^-|^2\Big)r\,{\mathrm d}r.
\end{align}

\medskip
Now we prove Lemma \ref{1.49} by contradiction.
Suppose that there exists a sequence of functions $\tilde {\psi}_l=(\tilde {\psi}_l^+, \tilde {\psi}_l^-)$
with $\|\tilde {\psi}_l\|_{\mathcal H_*}=1, \, l\in \mathbb{N}^+$ such that
$$
\langle\tilde {\psi}_l, Z_0\rangle_* =0,\, \forall\, l\in \mathbb{N}^+
\quad\mbox{and}\quad
\mathscr{B}\big({\tilde\psi}_l, {\tilde\psi}_l\big)\rightarrow 0
\ \mbox{as }  l\rightarrow +\infty.
$$
Let $\hat \psi$ be the weak limit of ${\tilde\psi}_l$ in the sense of $\|\cdot\|_{\mathcal H_*}$.
We claim that $\hat \psi\neq 0$. Indeed, ${\tilde\psi}_l\rightarrow\hat \psi$ locally in $L^2$ sense by compactly embedding theorem.
Hence, if $\hat \psi\equiv0$, we would have
$$
\int_0^R\big(|{\tilde\psi}_l^+|^2+|{\tilde\psi}_l^-|^2\big)r\,{\mathrm d}r\rightarrow 0.
$$
The estimate (\ref{1.32}) implies that $\|\tilde\psi_l\|_{\mathcal H_*}\rightarrow 0$, which is impossible, thus $\hat \psi\equiv0 $ does not hold.
Strong $L^2$ convergence over compacts and weak semi-continuity  of $L^2$-norm imply that
$$
\mathscr{B}(\hat\psi,\hat\psi)=0.
$$
Then we have $\hat\psi=CZ_0$.  Since $\{\tilde {\psi}_l\}_{l=1}^\infty $ are uniformly bounded in $\|\cdot\|_{\mathcal H_*}$ norm, then we conclude the weak convergence in  $\|\cdot\|_{\mathcal H_*}$ norm, which  implies that
 $$\langle\tilde {\psi}_l, Z_0\rangle_*\,\rightarrow\,  \langle\hat\psi, Z_0\rangle_*=0\quad \mbox{as }  l\rightarrow +\infty,$$
 so $C=0$. That is a contradiction, so we prove the lemma.\qed

\medskip
Define the closed subspace $\mathcal H_0$ of $ \mathcal H_*$  in the form
\begin{align*}
\mathcal H_0=\{\tilde\psi\in  \mathcal H_*\,:\, \langle \tilde\psi,Z_0\rangle_*=0\}.
\end{align*}
Using Lemma \ref{1.49}, we can know that the functional  $J_1$ is continuous, coercivity, and convex.
Thus there exists  a minimizer $\tilde\psi$ of the energy functional  $J_1$ in $\mathcal H_0$.
Indeed, by using orthogonal projection onto  the closed subspace $\mathcal H_0$, together with the  orthogonal  condition  (\ref{a.1}), we can know that $\tilde\psi$ is also the minimizer  of  $J_1$ in $\mathcal H_*$.
We omit the details here and
the readers can refer  to  the proof of Theorem 1.2 in \cite{DFK}.
It is straightforward to check that the inherited solution $\psi_1^1$ of \,${\mathcal L}(\psi_1^1)=h_1^1$\,\,indeed satisfies
\begin{align}\label{1.52}
\|\psi_1^1\|_{\mathcal H}^2\leq C\int_{\mathbb{R}^2}|h|^2(1+r^{2+\sigma}).
\end{align}

\medskip
In a similar way, we can find a solution $\psi_1^2$ to the equation ${\mathcal L}(\psi_1^2)=h_1^2$ with analogous estimate.
This solves  the equation (\ref{1.31}) for $j=1$.

\medskip
\noindent{\textbf{  The mode $j\geq  2$.}}  For the case $j\geq 2$, we define the closed subspace $\mathcal H^\bot$  of all functions $\psi^\bot\in  \mathcal H$ that can be
written in the form of
\begin{align*}
\psi^\bot=\sum_{j\geq 2}\psi_j^1+\sum_{j\geq 2}\psi_j^2,
\end{align*}
where $\psi_j^1,\, \psi_j^2, \forall\, j\geq 2$ are given in (\ref{decompositionphi}).
For a given function $$
h^\bot=\sum_{j\geq 2}h_j^1+\sum_{j\geq 2}h_j^2.
$$
where $h_j^1,\, h_j^2, \forall\, j\geq 2$ are given in (\ref{decompositionh}),
we now consider the equation
\begin{align}\label{bot}
{\mathcal L}(\psi^\bot)=h^\bot.
\end{align}
 In order    to solve the equation \eqref{bot}, we need to prove   the coercivity of the quadratic form
\begin{align}\label{quadratic}
\boldsymbol{B}(\psi^\bot,\psi^\bot)\geq C\|\psi^\bot\|_{\mathcal {H}}^2,
\end{align}
for some constant $C>0$.
We will adopt  relevant technique from the work  \cite{DFK} to prove \eqref{quadratic}.

 Recall the decompositions in (\ref{decompositionphi})-(\ref{decompositionh0}),
by  accurate  calculations of all terms in $\boldsymbol{B}(\psi,\psi)$,  we can get a similar result as (1.16) in \cite{DFK}, i.e.
 \begin{align}\label{1.28}
\boldsymbol{B}(\psi,\psi)=  \boldsymbol{B} (\psi_0, \psi_0)+ \sum_{j=1}^{\infty}\boldsymbol{B}(\psi_j^1, \psi_j^1)+\sum_{j=1}^{\infty}\boldsymbol{B}(\psi_j^2, \psi_j^2).
 \end{align}

\medskip
Next, for later use, we make efforts to translate  the  quadratic form $\boldsymbol{B}(\psi_j^\ell, \psi_j^\ell)$  to an equivalent quadratic form
$\mathcal{B}_j^\ell(\varphi,\varphi)$ in (\ref{1.37}),  which is defined on a space involving radial  functions of real vector values.
This will be fulfilled in two steps.
\medskip
\\
\noindent{\bf {Step 1.}}
We  define $(\varphi^+,\varphi^-)$  by the relation
\begin{align}\label{varphi}
\psi=(\psi^+,\psi^-)=(iw^+\varphi^+,iw^-\varphi^-),
\end{align}
where $w=(w^+,w^-)$ is the radially symmetric vortex solution  with degree pair $(n_+,n_-)=(1,1)$, see (\ref{(w^+,w^-)}). Then we introduce the quadratic form for $(\varphi^+,\varphi^-)$
\begin{align}\label{M}
{\mathbb M}(\varphi,\varphi)&=\boldsymbol{B}(iw^+\varphi^+,iw^-\varphi^-)\nonumber
\\[1mm]
&=\int_{\mathbb R^2}\Big[\, {U^+}^2|\nabla{\varphi^+}|^2+{U^-}^2|\nabla{\varphi^-}|^2\,\Big]
\,\,
-2{\mathbf {Re}}\int_{\mathbb R^2}\frac{1}{r^2}\Big[\,i{U^+}^2\frac{\partial\varphi^+}{\partial\theta}+i{U^-}^2
\frac{\partial\varphi^-}{\partial\theta}\,\Big]\nonumber
\\[1mm]
&~~~~\,+\,\int_{\mathbb R^2}\Big[\,2A_+{U^+}^4|\varphi_2^+|^2+2A_-{U^-}^4|\varphi_2^-|^2\,\Big]
\,+\,
4B\int_{\mathbb R^2}{U^+}^2{U^-}^2\varphi_2^+\varphi_2^-,
\end{align}
where we have used the convention
$$
\varphi=(\varphi^+,\varphi^-)=(\varphi_1^++i\varphi_2^+,\varphi_1^-+i\varphi_2^-).
$$
In fact,  the result in  (\ref{M})  is similar to the result in \cite{DFK} which in the case  the quadratic form are defined for complex-valued scalar functions.
We here omit the details of the proof for concise.
We now make a decomposition of ${\mathbb M}(\varphi,\varphi)$.
By the relations
\begin{align} \label{psi0}
\psi_0=(iw^+{\varphi_0}^+,iw^-{\varphi_0}^-),\quad \psi_j^\ell=(iw^+{\varphi_j^\ell}^+,iw^-{\varphi_j^\ell}^-),~\ell=1,2, \, j\in \mathbb{N}^+,
\end{align}
where $\psi_0, \psi_j^\ell, \ell=1,2,\, j\in \mathbb{N}^+,$ are the functions in \eqref{phi0}-\eqref{decompositionh0},
naturally, we have
\begin{align*}
\varphi=\varphi_0+\sum_{j=1}^\infty{\varphi_j^1}+\sum_{j=1}^\infty{\varphi_j^2},
\end{align*}
where
$$
\varphi_0=({\varphi_0}^+,{\varphi_0}^-),~~~\varphi_j^\ell=({\varphi_j^\ell}^+,{\varphi_j^\ell}^-),~~~\ell=1,2,\, j\in \mathbb{N}^+.
$$
Using (\ref{1.28}), the decomposition is
\begin{align*}
{\mathbb M}(\varphi,\varphi)=  {\mathbb M} (\varphi_0,\varphi_0)+ \sum_{j=1}^{\infty}{\mathbb M}(\varphi_j^1,\varphi_j^1)+\sum_{j=1}^{\infty}{\mathbb M}(\varphi_j^2,\varphi_j^2).
 \end{align*}

\noindent{\bf{Step 2.}}
To proceed, let us set
\begin{align}
 &{\varphi_0}^{\pm}={\varphi_{01}}^{\pm}+i{\varphi_{02}}^{\pm},\nonumber
\\[1mm]
  &{\varphi_j^1}^{\pm}={\varphi_{j1}^{1}}^{\pm}\cos{j\theta}+i{\varphi_{j2}^{1}}^{\pm}\sin{j\theta},
\\[1mm]
  &{\varphi_j^2}^{\pm}={\varphi_{j1}^2}^{\pm}\sin{j\theta}+i{\varphi_{j2}^2}^{\pm}\cos{j\theta},
\nonumber
  \end{align}
where $$
\varphi_0=({\varphi_0}^+,{\varphi_0}^-),~~~\varphi_j^\ell=({\varphi_j^\ell}^+,{\varphi_j^\ell}^-),~~~\ell=1,2,\, j\in \mathbb{N}^+.
$$
 are defined in \eqref{psi0}.
We then define the real vectors with
\begin{align}\label{mathbfV}
{\mathbf V}_{j}^\ell=({{\mathbf V}_j^\ell}^+,{{\mathbf V}_j^\ell}^-)=\Big(\big({\varphi_{j1}^{\ell}}^{+}, {\varphi_{j2}^{\ell}}^{+}\big),\, \big({\varphi_{j1}^{\ell}}^{-}, {\varphi_{j2}^{\ell}}^{-}\big)\Big).
\end{align}
According to \eqref{M},  we consider the quadratic forms $ \mathcal{B}_j^\ell({\mathbf V},{\mathbf V})$  for $j\in\mathbb{N}^+$ and $\ell=1,2$,
 \begin{align}\label{1.48}
 \mathcal{B}_j^\ell({\mathbf V},{\mathbf V})=&\int_0^\infty \Big[\,{U^+}^2|{{\mathbf V}^+}'|^2+{U^-}^2|{{\mathbf V}^-}'|^2+{U^+}^2{M_j^\ell}{\mathbf V}^+\cdot{\mathbf V}^++{U^-}^2{M_j^\ell}{\mathbf V}^-\cdot{\mathbf V}^-\,\Big]r\,{\mathrm d}r
 \nonumber\\
 &+
\int_0^\infty \Big[\,2A_+{U^+}^4{{\mathbf V}^+_2}^2+ 2A_-{U^-}^4{{\mathbf V}^-_2}^2+4B{U^+}^2{U^-}^2{\mathbf V}^+_2{\mathbf V}^-_2\,\Big]r\,{\mathrm d}r,
 \end{align}
 where the vector ${\mathbf V}=({\mathbf V}^+,{\mathbf V}^-)=\big(({\mathbf V}_1^+, {\mathbf V}_2^+), ({\mathbf V}_1^-, {\mathbf V}_2^-)\big):[0,\infty]\rightarrow \mathbb{R}^2\times{\mathbb{R}^2}$
 and
\begin{equation}\label{1.39}
{M_j^\ell}=\frac{1}{r^2}{
\left[ \begin{array}{ccc}
j^2& (-1)^{\ell+1}2j
\\[1mm]
(-1)^{\ell+1}2j & j^2
\\
\end{array}
\right ]}.
\end{equation}
It is easy to check that
\begin{align}\label{1.37}
\boldsymbol{B}(\psi_j^\ell, \psi_j^\ell)={\mathbb M}(\varphi_j^\ell,\varphi_j^\ell)=\pi \mathcal{B}_j^\ell({\mathbf V}_{j}^\ell,{\mathbf V}_{j}^\ell).
\end{align}
For $j\geq 2$, $\ell=1, 2$,  matrix
  ${M_j^\ell}$  given in (\ref{1.39})  satisfy
 \begin{equation*}
 (M_j^\ell- M_1^\ell){\mathbf V}^\pm\cdot{\mathbf V}^\pm=
\frac{j-1}{r^2}{
\left[ \begin{array}{ccc}
{j+1}& (-1)^{\ell+1}2
\\[1mm]
(-1)^{\ell+1}2 & {j+1} \\
\end{array}
\right ]}{\mathbf V}^\pm\cdot{\mathbf V}^\pm\geq \frac{(j-1)^2}{r^2} \, |{\mathbf V}^\pm|^2.
\end{equation*}
Therefore \begin{align}\label{1.36}
 \mathcal{B}_j^\ell({\mathbf V}, {\mathbf V})\geq
\int_0^\infty\frac{(j-1)^2}{r^2}\,\Big( {U^+}^2\,  |{\mathbf V}^+|^2+{U^-}^2\,  |{\mathbf V}^-|^2\Big)\, r\,{\mathrm d}r.
\end{align}
Together with  (\ref{1.28}), (\ref{1.37}), (\ref{1.36}), we can get
\begin{align*}
\boldsymbol{B}(\psi^\bot,\psi^\bot)&= \sum_{j\geq2}^{\infty}\boldsymbol{B}(\psi_j^1,\psi_j^1)+\sum_{j\geq 2}^{\infty}\boldsymbol{B}(\psi_j^2,\psi_j^2)
\\[1mm]
&=\sum_{j\geq2}^{\infty}{\mathbb M}(\varphi_j^1,\varphi_j^1)+\sum_{j\geq 2}^{\infty}{\mathbb M}(\varphi_j^2,\varphi_j^2)
\\[1mm]
&=\pi\, \sum_{j\geq2}^{\infty} \mathcal{B}_j^1({\mathbf V}_j^1,{\mathbf V}_j^1)+\,\pi \sum_{j\geq 2}^{\infty} \mathcal{B}_j^2({\mathbf V}_j^2,{\mathbf V}_j^2)
\\[1mm]
&\geq\, \pi\, \sum_{\ell=1}^2\sum_{j\geq2}^{\infty} \int_0^\infty \frac{(j-1)^2}{r^2}\Big[\,\big|{{\mathbf V}_j^\ell}^+\big|^2{U^+}^2+\big|{{\mathbf V}_j^\ell}^-\big|^2{U^-}^2 \,\Big]
\\[1mm]
&\geq C\int_{\mathbb{R}^2}\frac{|\psi^\bot|^2}{r^2}.
 \end{align*}
 where $ \psi_j^\ell, \varphi_j^\ell, {\mathbf V}_j^\ell$ are defined in \eqref{phi0},\, \eqref{decompositionh0},\, \eqref{psi0},\,\eqref{mathbfV}.
Then we can get the coercivity of the quadratic form $\boldsymbol{B}(\psi^\bot,\psi^\bot)$, see \eqref{quadratic}.

Since the functional  is continuous, coercive, and strictly convex in $\mathcal {H}^\bot$,
it is easy to conclude that for $h=h^\bot$ in (\ref{functionalJ}) the functional $J(\psi)$ has a
 minimizer $\psi^\bot$ in $\mathcal {H}^\bot$.

\medskip
Similar as the estimate (\ref{1.52}), we can get
\begin{align}\label{1.55}
\|\psi^\bot\|_{ \mathcal H}^2\leq C\int_{\mathbb{R}^2}|h|^2(1+r^{2+\sigma}).
\end{align}

\medskip
\noindent{\textbf{ The mode $j= 0$.}}
In this part, we will solve the equation \eqref{1.47}.
Setting $ \chi=(\chi^+,\chi^-)$, ${\mathfrak h}=({\mathfrak h}^+,{\mathfrak h}^-)$ by the relations
$$
\phi_0^\pm=iw^\pm\chi^\pm=iw^\pm(\chi^\pm_1+i\chi^\pm_2),\qquad
h_0^\pm=iw^\pm {\mathfrak h}^\pm=iw^\pm({\mathfrak h}_1^\pm+i{\mathfrak h}_2^\pm),
$$
we  can get an ODE system from the equation ${\mathcal L}(\phi_0)=h_0$,
\begin{align}
&{\chi^+_1}''+\Big(\frac{2{U^+}'}{U^+}+\frac{1}{r}\Big){\chi^+_1}'={\mathfrak h}_1^+\label{1.29},
\\[1mm]
&{\chi^-_1}''+\Big(\frac{2{U^-}'}{U^-}+\frac{1}{r}\Big){\chi^-_1}'={\mathfrak h}_1^-\label{1.30},
\\[1mm]
&{\chi^+_2}''+\Big(\frac{2{U^+}'}{U^+}+\frac{1}{r}\Big){\chi^+_2}'-2A_+{U^+}^2\chi^+_2-2B{U^-}^2\chi^-_2={\mathfrak h}_2^+,\label{1.50}
\\[1mm]
&{\chi^-_2}''+\Big(\frac{2{U^-}'}{U^-}+\frac{1}{r}\Big){\chi^-_2}'-2A_-{U^-}^2\chi^-_2-2B{U^+}^2\chi^+_2={\mathfrak h}_2^-.\label{1.51}
\end{align}
We can get the properties by using the assumptions (\ref{1.35}) and (\ref{03}) made on $h$
\begin{align}
\int_0^{\infty}r{U^{\pm}}^2{\mathfrak h}_1^{\pm}=0,
\qquad
\int_0^{\infty}r{U^{\pm}}^2|{\mathfrak h}_1^\pm|^2(1+r^{2+\sigma}) \,{\mathrm d}r\leq C\int_{{\mathbb R}^2}|h|^2(1+r^{2+\sigma}).
\end{align}
Using the variation of parameters,~we get the solutions of   (\ref{1.29})-(\ref{1.30})
\begin{align*}
\chi^+_1=-\int_r^\infty\frac{1}{s{U^+(s)}^2}\mathrm ds\int_0^s{\mathfrak h}_1^+(t){U^+(t)}^2\,\mathrm dt,
\\[1mm]
\chi^-_1=-\int_r^\infty\frac{1}{s{U^-(s)}^2}\mathrm ds\int_0^s{\mathfrak h}_1^-(t){U^-(t)}^2\,\mathrm dt,
\end{align*}
and the estimate
\begin{align*}
\int_0^\infty\Big(|{\chi^+_1}'(r)|^2+|{\chi^-_1}'(r)|^2\Big)r\,{\mathrm d}r\leq C\int _{\mathbb{R}^2}|h|^2(1+r^{2+\sigma}).
\end{align*}
On the other hand, we can solve equations (\ref{1.50}) and (\ref{1.51}) by finding a minimizer of the functional
\begin{align*}
J_2(\chi^+_2,\chi^-_2)
=&\,\frac{1}{2}\int_0^\infty\Big[\,\big|{\chi^+_2}'(r)\big|^2{U^+}^2+\big|{\chi^-_2}'(r)\big|^2{U^-}^2+2A_+{U^+}^4{\chi^+_2}^2+2A_-{U^-}^4{\chi^-_2}^2\,\Big]r\,{\mathrm d}r
\nonumber
\\[1mm]
&+\int_0^\infty
\Big[\,4B{U^+}^2{U^-}^2\chi^+_2\chi^-_2+{\mathfrak h}_2^+\chi^+_2{U^+}^2+{\mathfrak h}_2^-\chi^-_2{U^-}^2\,\Big]r\,{\mathrm d}r.
\end{align*}
With these works, we  obtain the estimate for $\phi_0$
\begin{align}\label{1.54}
\|\psi_0\|_{\mathcal H}^2\leq C\int_{\mathbb{R}^2}|h|^2(1+r^{2+\sigma}).
\end{align}

\medskip
As a conclusion, together with Steps 1-3 and the estimates (\ref{1.52}), (\ref{1.55}), (\ref{1.54}), we can obtain a solution $\psi$ of (\ref{1.22}) with the required properties. The fact that the solutions of (\ref{1.22}) can be written  as in (\ref{1.56}) is a direct corollary of Theorem \ref{theorem1.1}. This complete the proof of Theorem \ref{theorem1.2}.
\qed

\bigskip
\bigskip
%\begin{acknowledgements}
{\bf Acknowledgements: }
Jun Yang is supported by  the grants CCNU18CXTD04  and NSFC (No. 11771167).
\qed
%\end{acknowledgements}

%\begin{thebibliography}{99}

 \end{document}